\documentclass[12pt]{article}

\newcommand{\be}{\begin{equation}}
\newcommand{\ee}{\end{equation}}
\newcommand{\bea}{\begin{eqnarray}}
\newcommand{\eea}{\end{eqnarray}}
\newcommand{\beaa}{\begin{eqnarray*}}
\newcommand{\eeaa}{\end{eqnarray*}}

\newcommand{\spc}{\hspace{0.2in}}
\newcommand{\g}{{\frak g}}

\newcommand{\C}{{\mathbb C}}
\newcommand{\I}{1}
\newcommand{\transp}{{\bf P}}

\newcommand{\hsp}{\hspace{0.2in}}

\advance\textheight by \topskip
\usepackage{amssymb}
\usepackage{pstricks}

\pslinewidth .8pt \psset{unit=1.5cm} \psset{fillcolor=white}

\newcommand{\Ru}{\begin{pspicture}(0,.2)(.6,.7)
\psline(0,0)(.14,.14) \psline(.46,.46)(.6,.6)
\psline(0,.6)(.14,.46)\psline(.46,.14)(.6,0)\pscircle(.3,.3){.2}
\rput[bl](-.15,.53){$a$}\rput[bl](-.15,-.1){$b$}
\rput[bl](.6,-.1){$d$} \rput[bl](.6,.53){$c$} \rput(.3,.3){$u$}
\end{pspicture} }

\newcommand{\RuD}{\begin{pspicture}(0,.2)(.6,.7)
\psline(0,0)(.14,.14)\psline[arrowsize=2pt
3,arrowinset=0.75]{<-}(0,0)(.14,.14)
\psline(.46,.46)(.6,.6)\psline[arrowsize=2pt
3,arrowinset=0.75]{<-}(.46,.46)(.6,.6)
\psline(0,.6)(.14,.46)\psline[arrowsize=2pt
3,arrowinset=0.75]{<-}(0,.6)(.14,.46)
\psline(.46,.14)(.6,0)\psline[arrowsize=2pt
3,arrowinset=0.75]{<-}(.46,.14)(.6,0)\pscircle(.3,.3){.2}
\rput[bl](-.15,.53){$a$}\rput[bl](-.15,-.1){$b$}
\rput[bl](.6,-.1){$d$} \rput[bl](.6,.53){$c$} \rput(.3,.3){$u$}
\end{pspicture} }

\newcommand{\RuDc}{\begin{pspicture}(0,.2)(.6,.7)
\psline(0,0)(.14,.14)\psline[arrowsize=2pt
3,arrowinset=0.75]{<-}(0,0)(.14,.14)
\psline(.46,.46)(.6,.6)\psline[arrowsize=2pt
3,arrowinset=0.75]{<-}(.46,.46)(.6,.6)
\psline(0,.6)(.14,.46)\psline[arrowsize=2pt
3,arrowinset=0.75]{->}(0,.6)(.14,.46)
\psline(.46,.14)(.6,0)\psline[arrowsize=2pt
3,arrowinset=0.75]{->}(.46,.14)(.6,0)\pscircle(.3,.3){.2}
\rput[bl](-.15,.53){$a$}\rput[bl](-.15,-.1){$b$}
\rput[bl](.6,-.1){$d$} \rput[bl](.6,.53){$c$} \rput(.3,.3){$u$}
\end{pspicture} }

\newcommand{\RUu}
{\begin{pspicture}(0,.2)(.6,.7) \psline(0,0)(.6,0)
\psline(0,.3)(.13,.34) \psline(.47,.56)(.6,.6)
\psline(0,.6)(.13,.56)\psline(.47,.34)(.6,.3)
\pscircle(.3,.45){.2} \rput(.3,.45){$u$}
\end{pspicture} }

\newcommand{\RUv}
{\begin{pspicture}(0,.2)(.6,.7) \psline(0,0)(.6,0)
\psline(0,.3)(.13,.34) \psline(.47,.56)(.6,.6)
\psline(0,.6)(.13,.56)\psline(.47,.34)(.6,.3)
\pscircle(.3,.45){.2} \rput(.3,.45){$v$}
\end{pspicture} }

\newcommand{\RUuv}
{\begin{pspicture}(0,.2)(.6,.7) \psline(0,0)(.6,0)
\psline(0,.3)(.13,.34) \psline(.47,.56)(.6,.6)
\psline(0,.6)(.13,.56)\psline(.47,.34)(.6,.3)
\pscircle(.3,.45){.2} \rput(.3,.45){${\scriptstyle u\!+\!v}$}
\end{pspicture} }

\newcommand{\RDu}{\begin{pspicture}(0,.2)(.6,.7)
\psline(0,0)(.13,.04)\psline(.47,.26)(.6,.3)
\psline(0,.3)(.13,.26)\psline(.47,.04)(.6,0) \psline(0,.6)(.6,.6)
\pscircle(.3,.15){.2} \rput(.3,.15){$u$}
\end{pspicture} }

\newcommand{\RDv}{\begin{pspicture}(0,.2)(.6,.7)
\psline(0,0)(.13,.04)\psline(.47,.26)(.6,.3)
\psline(0,.3)(.13,.26)\psline(.47,.04)(.6,0) \psline(0,.6)(.6,.6)
\pscircle(.3,.15){.2} \rput(.3,.15){$v$}
\end{pspicture} }

\newcommand{\RDuv}{\begin{pspicture}(0,.2)(.6,.7)
\psline(0,0)(.13,.04)\psline(.47,.26)(.6,.3)
\psline(0,.3)(.13,.26)\psline(.47,.04)(.6,0) \psline(0,.6)(.6,.6)
\pscircle(.3,.15){.2} \rput(.3,.15){${\scriptstyle u\!+\!v}$}
\end{pspicture} }

\newcommand{\obj}
{\begin{pspicture}(0,.2)(.6,.7) \psline(0,0)(.2,.2)
\psline(.4,.4)(.6,.6) \psline(0,.6)(.2,.4)\psline(.4,.2)(.6,0)
\psline(0,.3)(.16,.3) \psline(.44,.3)(.6,.3)\pscircle(.3,.3){.14}
\psline[arrowsize=2pt 3,arrowinset=0.75]{<-}(0,0)(.2,.2)
\psline[arrowsize=2pt 3,arrowinset=0.75]{<-}(.4,.4)(.6,.6)
\psline[arrowsize=2pt
3,arrowinset=0.75]{<-}(0,.6)(.2,.4)\psline[arrowsize=2pt
3,arrowinset=0.75]{<-}(.4,.2)(.6,0) \psline[arrowsize=2pt
3,arrowinset=0.75]{<-}(0,.3)(.16,.3) \psline[arrowsize=2pt
3,arrowinset=0.75]{<-}(.44,.3)(.6,.3)
\end{pspicture} }

\newcommand{\objj}
{\begin{pspicture}(0,.2)(.6,.7) \psline(0,0)(.2,.2)
\psline(.4,.4)(.6,.6) \psline(0,.6)(.2,.4)\psline(.4,.2)(.6,0)
\psline(0,.3)(.16,.3)
\psline(.44,.3)(.6,.3)\pscircle(.3,.3){.14}\pscircle(.3,.3){.1}
\psline[arrowsize=2pt 3,arrowinset=0.75]{<-}(0,0)(.2,.2)
\psline[arrowsize=2pt 3,arrowinset=0.75]{<-}(.4,.4)(.6,.6)
\psline[arrowsize=2pt
3,arrowinset=0.75]{<-}(0,.6)(.2,.4)\psline[arrowsize=2pt
3,arrowinset=0.75]{<-}(.4,.2)(.6,0) \psline[arrowsize=2pt
3,arrowinset=0.75]{<-}(0,.3)(.16,.3) \psline[arrowsize=2pt
3,arrowinset=0.75]{<-}(.44,.3)(.6,.3)
\end{pspicture} }

\newcommand{\objh}
{\begin{pspicture}(0,.2)(.6,.7) \psline(0,0)(.2,.2)
\psline(.4,.4)(.6,.6) \psline(0,.6)(.2,.4)\psline(.4,.2)(.6,0)
\psline(0,.3)(.16,.3)
\psline(.44,.3)(.6,.3)\pscircle(.3,.3){.14}\pscircle*(.3,.3){.1}
\end{pspicture} }

\newcommand{\G}{\begin{pspicture}(0,.2)(.6,.6)
\psline[arrowsize=2pt 6,arrowinset=.4]{<-}(.15,.3)(.6,.3)
\psline[dotsize=2pt 5.2]{o-}(.36,.3)(.6,.3) \psline(0,.3)(.6,.3)
\end{pspicture} }

\newcommand{\GT}{\begin{pspicture}(0,.2)(.6,.6)
\psline[arrowsize=2pt 6,arrowinset=.4]{->}(0,.3)(.48,.3)
\psline[dotsize=2pt 5.2]{o-}(.27,.3)(.6,.3) \psline(0,.3)(.6,.3)
\end{pspicture} }

\newcommand{\GTbar}{\begin{pspicture}(0,.2)(.6,.6)
\psline[arrowsize=2pt 6,arrowinset=.4]{->}(0,.3)(.48,.3)
\psline[dotsize=2pt 5.2]{o-}(.27,.3)(.6,.3)
\psline(.17,.4)(.37,.4)\psline(0,.3)(.6,.3)
\end{pspicture} }

\newcommand{\GG}{\begin{pspicture}(0,.2)(.6,.6)
\psline[arrowsize=2pt 6,arrowinset=0.4]{<-}(.15,0)(.6,0)
\psline[arrowsize=2pt 6,arrowinset=0.4]{<-}(.15,.6)(.6,.6)
\psline[dotsize=2pt 5.2]{o-}(.36,0)(.6,0) \psline[dotsize=2pt
5.2]{o-}(.36,.6)(.6,.6) \psline(0,0)(.6,0) \psline(0,.6)(.6,.6)
\end{pspicture} }

\newcommand{\GGG}{\begin{pspicture}(0,.2)(.6,.6)
\psline[arrowsize=2pt 6,arrowinset=0.4]{<-}(.15,0)(.6,0)
\psline[arrowsize=2pt 6,arrowinset=0.4]{<-}(.15,.3)(.6,.3)
\psline[arrowsize=2pt 6,arrowinset=0.4]{<-}(.15,.6)(.6,.6)
\psline[dotsize=2pt 5.2]{o-}(.36,0)(.6,0)\psline[dotsize=2pt
5.2]{o-}(.36,.3)(.6,.3) \psline[dotsize=2pt
5.2]{o-}(.36,.6)(.6,.6) \psline(0,0)(.6,0)\psline(0,.3)(.6,.3)
\psline(0,.6)(.6,.6)
\end{pspicture} }

\newcommand{\GGGG}{\begin{pspicture}(0,.2)(.6,.6)
\psline[arrowsize=2pt 6,arrowinset=0.4]{<-}(.15,0)(.6,0)
\psline[arrowsize=2pt 6,arrowinset=0.4]{<-}(.15,.2)(.6,.2)
\psline[arrowsize=2pt 6,arrowinset=0.4]{<-}(.15,.4)(.6,.4)
\psline[arrowsize=2pt 6,arrowinset=0.4]{<-}(.15,.6)(.6,.6)
\psline[dotsize=2pt 5.2]{o-}(.36,0)(.6,0)\psline[dotsize=2pt
5.2]{o-}(.36,.2)(.6,.2)\psline[dotsize=2pt 5.2]{o-}(.36,.4)(.6,.4)
\psline[dotsize=2pt 5.2]{o-}(.36,.6)(.6,.6)
\psline(0,0)(.6,0)\psline(0,.2)(.6,.2)\psline(0,.4)(.6,.4)
\psline(0,.6)(.6,.6)
\end{pspicture} }

\newcommand{\suAA}{\begin{pspicture}(0,.2)(.6,.6)
\psline[arrowsize=2pt 2.5,arrowinset=0.75]{<-}(0,.3)(.6,.3)
\end{pspicture} }

\newcommand{\suBA}{\begin{pspicture}(0,.2)(.6,.6)
\psline[arrowsize=2pt 3,arrowinset=0.75]{<-}(0,0)(.6,0)
\psline[arrowsize=2pt 3,arrowinset=0.75]{<-}(0,.6)(.6,.6)
\psline(0,0)(.6,0) \psline(0,.6)(.6,.6)
\end{pspicture} }

\newcommand{\suBB}{\begin{pspicture}(0,.2)(.6,.6)
\psline[arrowsize=2pt 3,arrowinset=0.75]{<-}(0,0)(.6,.6)
\psline[arrowsize=2pt 3,arrowinset=0.75]{<-}(0,.6)(.6,0)
\psline(0,0)(.6,.6) \psline(0,.6)(.6,0)
\end{pspicture} }

\newcommand{\suBC}{\begin{pspicture}(0,.2)(.6,.6)
\psline[arrowsize=2pt 3,arrowinset=0.75]{<-}(0,0)(.6,.6)
\psline[arrowsize=2pt 3,arrowinset=0.75]{->}(0,.6)(.6,0)
\psline(0,0)(.6,.6) \psline(0,.6)(.6,0)
\end{pspicture} }

\newcommand{\suBD}{\begin{pspicture}(0,.2)(.8,.7)
\psarc(0,.3){.3}{270}{90} \psarc(.8,.3){.3}{90}{270}
\psarc[arrowsize=2pt 3,arrowinset=0.75]{<-}(0,.3){.3}{270}{90}
\psarc[arrowsize=2pt 3,arrowinset=0.75]{->}(.8,.3){.3}{90}{270}
\end{pspicture} }

\newcommand{\suBE}{\begin{pspicture}(0,.2)(.6,.6)
\psline[arrowsize=2pt 3,arrowinset=0.75]{->}(0,0)(.6,0)
\psline[arrowsize=2pt 3,arrowinset=0.75]{<-}(0,.6)(.6,.6)
\psline(0,0)(.6,0) \psline(0,.6)(.6,.6)
\end{pspicture} }

\newcommand{\suBF}{\begin{pspicture}(0,.2)(.8,.7)
\psarc(0,.3){.3}{270}{90} \psarc(.8,.3){.3}{90}{270}
\psarc[arrowsize=2pt 3,arrowinset=0.75]{->}(0,.3){.3}{270}{90}
\psarc[arrowsize=2pt 3,arrowinset=0.75]{->}(.8,.3){.3}{90}{270}
\end{pspicture} }

\newcommand{\suCA}{\begin{pspicture}(0,.2)(.5,.6)
\psline[arrowsize=2pt 2.5,arrowinset=0.75]{<-}(0,0)(.5,0)
\psline[arrowsize=2pt 2.5,arrowinset=0.75]{<-}(0,.3)(.5,.3)
\psline[arrowsize=2pt 2.5,arrowinset=0.75]{<-}(0,.6)(.5,.6)
\psline(0,0)(.5,0) \psline(0,.3)(.5,.3) \psline(0,.6)(.5,.6)
\end{pspicture} }

\newcommand{\suCB}{\begin{pspicture}(0,.2)(.5,.6)
\psline[arrowsize=2pt 2.5,arrowinset=0.75]{<-}(0,0)(.5,0)
\psline[arrowsize=2pt 2.5,arrowinset=0.75]{<-}(0,.3)(.5,.6)
\psline[arrowsize=2pt 2.5,arrowinset=0.75]{<-}(0,.6)(.5,.3)
\psline(0,0)(.5,0) \psline(0,.3)(.5,.6) \psline(0,.6)(.5,.3)
\end{pspicture} }

\newcommand{\suCC}{\begin{pspicture}(0,.2)(.5,.6)
\psline[arrowsize=2pt 2.5,arrowinset=0.75]{<-}(0,0)(.5,.3)
\psline[arrowsize=2pt 2.5,arrowinset=0.75]{<-}(0,.3)(.5,0)
\psline[arrowsize=2pt 2.5,arrowinset=0.75]{<-}(0,.6)(.5,.6)
\psline(0,0)(.5,.3) \psline(0,.3)(.5,0) \psline(0,.6)(.5,.6)
\end{pspicture} }

\newcommand{\suCD}{\begin{pspicture}(0,.2)(.5,.6)
\psline[arrowsize=2pt 2.5,arrowinset=0.75]{<-}(0,0)(.5,.6)
\psline[arrowsize=2pt 2.5,arrowinset=0.75]{<-}(0,.3)(.5,.3)
\psline[arrowsize=2pt 2.5,arrowinset=0.75]{<-}(0,.6)(.5,0)
\psline(0,0)(.5,.6) \psline(0,.3)(.5,.3) \psline(0,.6)(.5,0)
\end{pspicture} }

\newcommand{\suCE}{\begin{pspicture}(0,.2)(.5,.6)
\psline[arrowsize=2pt 2.5,arrowinset=0.75]{<-}(0,0)(.5,.3)
\psline[arrowsize=2pt 2.5,arrowinset=0.75]{<-}(0,.3)(.5,.6)
\psline[arrowsize=2pt 2.5,arrowinset=0.75]{<-}(0,.6)(.5,0)
\psline(0,0)(.5,.3) \psline(0,.3)(.5,.6) \psline(0,.6)(.5,0)
\end{pspicture} }

\newcommand{\suCF}{\begin{pspicture}(0,.2)(.5,.6)
\psline[arrowsize=2pt 2.5,arrowinset=0.75]{<-}(0,0)(.5,.6)
\psline[arrowsize=2pt 2.5,arrowinset=0.75]{<-}(0,.3)(.5,0)
\psline[arrowsize=2pt 2.5,arrowinset=0.75]{<-}(0,.6)(.5,.3)
\psline(0,0)(.5,.6) \psline(0,.3)(.5,0) \psline(0,.6)(.5,.3)
\end{pspicture} }

\newcommand{\soAA}{\begin{pspicture}(0,.2)(.5,.6)
\psline(0,.3)(.5,.3)
\end{pspicture} }

\newcommand{\soBA}{\begin{pspicture}(0,.2)(.6,.7) \psline(0,0)(.6,0)
\psline(0,.6)(.6,.6)
\end{pspicture} }

\newcommand{\soBB}{\begin{pspicture}(0,.2)(.6,.7)
\psline(0,0)(.6,.6) \psline(0,.6)(.6,0)
\end{pspicture} }

\newcommand{\soBC}{\begin{pspicture}(0,.2)(.8,.7)
\psarc(0,.3){.3}{270}{90} \psarc(.8,.3){.3}{90}{270}
\end{pspicture} }

\newcommand{\soCA}{\begin{pspicture}(0,.2)(.6,.7)
\psline(0,0)(.6,0) \psline(0,.3)(.6,.3) \psline(0,.6)(.6,.6)
\end{pspicture} }

\newcommand{\soCB}{\begin{pspicture}(0,.2)(.6,.7)
\psline(0,0)(.6,0) \psline(0,.3)(.6,.6) \psline(0,.6)(.6,.3)
\end{pspicture} }

\newcommand{\soCC}{\begin{pspicture}(0,.2)(.6,.7)
\psline(0,0)(.6,.3) \psline(0,.3)(.6,0) \psline(0,.6)(.6,.6)
\end{pspicture} }

\newcommand{\soCD}{\begin{pspicture}(0,.2)(.6,.7)
\psline(0,0)(.6,.6) \psline(0,.3)(.6,.3) \psline(0,.6)(.6,0)
\end{pspicture} }

\newcommand{\soCE}{\begin{pspicture}(0,.2)(.6,.7)
\psline(0,0)(.6,.3) \psline(0,.3)(.6,.6) \psline(0,.6)(.6,0)
\end{pspicture} }

\newcommand{\soCF}{\begin{pspicture}(0,.2)(.6,.7)
\psline(0,0)(.6,.6) \psline(0,.3)(.6,0) \psline(0,.6)(.6,.3)
\end{pspicture} }

\newcommand{\soCG}{\begin{pspicture}(0,.2)(.6,.7)
\psline(0,0)(.6,0) \psarc(0,.45){.15}{270}{90}
\psarc(.6,.45){.15}{90}{270}
\end{pspicture} }

\newcommand{\soCH}{\begin{pspicture}(0,.2)(.6,.7)
\psline(0,.6)(.6,.6) \psarc(0,.15){.15}{270}{90}
\psarc(.6,.15){.15}{90}{270}
\end{pspicture} }

\newcommand{\soCI}{\begin{pspicture}(0,.2)(.6,.7)
\psline(0,.6)(.6,0) \psarc(0,.15){.15}{270}{90}
\psarc(.6,.45){.15}{90}{270}
\end{pspicture} }

\newcommand{\soCJ}{\begin{pspicture}(0,.2)(.6,.7)
\psline(0,0)(.6,.6) \psarc(0,.45){.15}{270}{90}
\psarc(.6,.15){.15}{90}{270}
\end{pspicture} }

\newcommand{\soCK}{\begin{pspicture}(0,.2)(.7,.7)
\psline(0,.3)(.7,.3) \psarc(0,.3){.3}{270}{90}
\psarc(.7,.3){.3}{90}{270}
\end{pspicture} }

\newcommand{\soCL}{\begin{pspicture}(0,.2)(.6,.7)
\psline(0,0)(.6,.3) \psarc(0,.45){.15}{270}{90}
\psarc(.6,.3){.3}{90}{270}
\end{pspicture} }

\newcommand{\soCM}{\begin{pspicture}(0,.2)(.6,.7)
\psline(0,.6)(.6,.3) \psarc(0,.15){.15}{270}{90}
\psarc(.6,.3){.3}{90}{270}
\end{pspicture} }

\newcommand{\soCN}{\begin{pspicture}(0,.2)(.6,.7)
\psline(0,.3)(.6,.6) \psarc(0,.3){.3}{270}{90}
\psarc(.6,.15){.15}{90}{270}
\end{pspicture} }

\newcommand{\soCO}{\begin{pspicture}(0,.2)(.6,.7)
\psline(0,.3)(.6,0) \psarc(0,.3){.3}{270}{90}
\psarc(.6,.45){.15}{90}{270}
\end{pspicture} }

\newcommand{\spAA}{\begin{pspicture}(0,.2)(.6,.6)
\psline[arrowsize=2pt 4,arrowinset=0]{<-}(.2,.3)(.6,.3)
\psline(0,.3)(.6,.3)
\end{pspicture} }

\newcommand{\spAAT}{\begin{pspicture}(0,.2)(.6,.6)
\psline[arrowsize=2pt 4,arrowinset=0]{->}(0,.3)(.4,.3)
\psline(0,.3)(.6,.3)
\end{pspicture} }

\newcommand{\spBA}{\begin{pspicture}(0,.2)(.6,.7) \psline(0,0)(.6,0)
\psline(0,.6)(.6,.6)
\end{pspicture} }

\newcommand{\spBB}{\begin{pspicture}(0,.2)(.6,.7)
\psline(0,0)(.6,.6) \psline(0,.6)(.6,0)
\end{pspicture} }

\newcommand{\spBC}{\begin{pspicture}(0,.2)(.8,.7)
\psarc[arrowsize=2pt 4,arrowinset=0]{->}(0,.3){.3}{270}{20}
\psarc(0,.3){.3}{20}{90} \psarc[arrowsize=2pt
4,arrowinset=0]{->}(.8,.3){.3}{90}{200}\psarc(.8,.3){.3}{200}{270}
\end{pspicture} }

\newcommand{\spCA}{\begin{pspicture}(0,.2)(.6,.7)
\psline(0,0)(.6,0) \psline(0,.3)(.6,.3) \psline(0,.6)(.6,.6)
\end{pspicture} }

\newcommand{\spCB}{\begin{pspicture}(0,.2)(.6,.7)
\psline(0,0)(.6,0) \psline(0,.3)(.6,.6) \psline(0,.6)(.6,.3)
\end{pspicture} }

\newcommand{\spCC}{\begin{pspicture}(0,.2)(.6,.7)
\psline(0,0)(.6,.3) \psline(0,.3)(.6,0) \psline(0,.6)(.6,.6)
\end{pspicture} }

\newcommand{\spCD}{\begin{pspicture}(0,.2)(.6,.7)
\psline(0,0)(.6,.6) \psline(0,.3)(.6,.3) \psline(0,.6)(.6,0)
\end{pspicture} }

\newcommand{\spCE}{\begin{pspicture}(0,.2)(.6,.7)
\psline(0,0)(.6,.3) \psline(0,.3)(.6,.6) \psline(0,.6)(.6,0)
\end{pspicture} }

\newcommand{\spCF}{\begin{pspicture}(0,.2)(.6,.7)
\psline(0,0)(.6,.6) \psline(0,.3)(.6,0) \psline(0,.6)(.6,.3)
\end{pspicture} }

\newcommand{\spCG}{\begin{pspicture}(0,.2)(.6,.7)
\psline(0,0)(.6,0) \psarc[arrowsize=2pt
3,arrowinset=0]{->}(0,.45){.15}{270}{30}\psarc(0,.45){.15}{30}{90}
\psarc[arrowsize=2pt
3,arrowinset=0]{->}(.6,.45){.15}{90}{210}\psarc(.6,.45){.15}{210}{270}
\end{pspicture} }

\newcommand{\spCH}{\begin{pspicture}(0,.2)(.6,.7)
\psline(0,.6)(.6,.6) \psarc[arrowsize=2pt
3,arrowinset=0]{->}(0,.15){.15}{270}{30}\psarc(0,.15){.15}{25}{90}
\psarc[arrowsize=2pt
3,arrowinset=0]{->}(.6,.15){.15}{90}{210}\psarc(.6,.15){.15}{205}{270}
\end{pspicture} }

\newcommand{\spCI}{\begin{pspicture}(0,.2)(.6,.7)
\psline(0,.6)(.6,0) \psarc[arrowsize=2pt
3,arrowinset=0]{->}(0,.15){.15}{270}{30}\psarc(0,.15){.15}{25}{90}
\psarc[arrowsize=2pt
3,arrowinset=0]{->}(.6,.45){.15}{90}{210}\psarc(.6,.45){.15}{205}{270}
\end{pspicture} }

\newcommand{\spCJ}{\begin{pspicture}(0,.2)(.6,.7)
\psline(0,0)(.6,.6) \psarc[arrowsize=2pt
3,arrowinset=0]{->}(0,.45){.15}{270}{30}\psarc(0,.45){.15}{25}{90}
\psarc[arrowsize=2pt
3,arrowinset=0]{->}(.6,.15){.15}{90}{210}\psarc(.6,.15){.15}{205}{270}
\end{pspicture} }

\newcommand{\spCK}{\begin{pspicture}(0,.2)(.7,.7)
\psline(0,.3)(.7,.3) \psarc[arrowsize=2pt
3,arrowinset=0]{->}(0,.3){.3}{270}{0}\psarc(0,.3){.3}{0}{90}
\psarc[arrowsize=2pt
3,arrowinset=0]{->}(.7,.3){.3}{90}{180}\psarc(.7,.3){.3}{180}{270}
\end{pspicture} }

\newcommand{\spCL}{\begin{pspicture}(0,.2)(.6,.7)
\psline(0,0)(.6,.3) \psarc[arrowsize=2pt
3,arrowinset=0]{->}(0,.45){.15}{270}{30}\psarc(0,.45){.15}{25}{90}
\psarc[arrowsize=2pt
3,arrowinset=0]{->}(.6,.3){.3}{90}{195}\psarc(.6,.3){.3}{195}{270}
\end{pspicture} }

\newcommand{\spCM}{\begin{pspicture}(0,.2)(.6,.7)
\psline(0,.6)(.6,.3) \psarc[arrowsize=2pt
3,arrowinset=0]{->}(0,.15){.15}{270}{30}\psarc(0,.15){.15}{25}{90}
\psarc[arrowsize=2pt 3,
arrowinset=0]{->}(.6,.3){.3}{90}{210}\psarc(.6,.3){.3}{205}{270}
\end{pspicture} }

\newcommand{\spCN}{\begin{pspicture}(0,.2)(.6,.7)
\psline(0,.3)(.6,.6) \psarc[arrowsize=2pt 3,
arrowinset=0]{->}(0,.3){.3}{270}{15}\psarc(0,.3){.3}{15}{90}
\psarc[arrowsize=2pt 3,
arrowinset=0]{->}(.6,.15){.15}{90}{210}\psarc(.6,.15){.15}{205}{270}
\end{pspicture} }

\newcommand{\spCO}{\begin{pspicture}(0,.2)(.6,.7)
\psline(0,.3)(.6,0) \psarc[arrowsize=2pt 3,
arrowinset=0]{->}(0,.3){.3}{270}{20}\psarc(0,.3){.3}{20}{90}
\psarc[arrowsize=2pt 3,
arrowinset=0]{->}(.6,.45){.15}{90}{210}\psarc(.6,.45){.15}{205}{270}
\end{pspicture} }

\newcommand{\SPCAA}{\begin{pspicture}(0,.1)(.4,.5)
\psline(0,0)(.4,0) \psline(0,.2)(.4,.2) \psline(0,.4)(.4,.4)
\psline[arrowsize=2pt 3,arrowinset=0]{<-}(.1,0)(.4,0)
\psline[arrowsize=2pt 3,arrowinset=0]{<-}(.1,.2)(.4,.2)
\psline[arrowsize=2pt 3,arrowinset=0]{<-}(.1,.4)(.4,.4)
\end{pspicture} }

\newcommand{\SPCB}{\begin{pspicture}(0,.2)(.6,.7)
\psline(0,0)(.6,0) \psline(0,.3)(.6,.6) \psline(0,.6)(.6,.3)
\psline[arrowsize=2pt 3,arrowinset=0]{<-}(.2,0)(.6,0)
\psline[arrowsize=2pt 3,arrowinset=0]{<-}(.05,.325)(.6,.6)
\psline[arrowsize=2pt 3,arrowinset=0]{<-}(.05,.575)(.6,.3)
\end{pspicture} }

\newcommand{\SPCC}{\begin{pspicture}(0,.2)(.6,.7)
\psline(0,0)(.6,.3) \psline(0,.3)(.6,0) \psline(0,.6)(.6,.6)
\psline[arrowsize=2pt 3,arrowinset=0]{<-}(.05,.025)(.6,.3)
\psline[arrowsize=2pt 3,arrowinset=0]{<-}(.05,.275)(.6,0)
\psline[arrowsize=2pt 3,arrowinset=0]{<-}(.2,.6)(.6,.6)
\end{pspicture} }

\newcommand{\SPCG}{\begin{pspicture}(0,.2)(.6,.7)
\psline(0,0)(.6,0) \psline[arrowsize=2pt
3,arrowinset=0]{<-}(.2,0)(.6,0)\psarc[arrowsize=2pt
3,arrowinset=0]{->}(0,.45){.15}{270}{30}\psarc(0,.45){.15}{30}{90}
\psarc[arrowsize=2pt
3,arrowinset=0]{->}(.6,.45){.15}{90}{210}\psarc(.6,.45){.15}{210}{270}
\end{pspicture} }

\newcommand{\SPCH}{\begin{pspicture}(0,.2)(.6,.7)
\psline(0,.6)(.6,.6) \psline[arrowsize=2pt
3,arrowinset=0]{<-}(.2,.6)(.6,.6) \psarc[arrowsize=2pt
3,arrowinset=0]{->}(0,.15){.15}{270}{30}\psarc(0,.15){.15}{25}{90}
\psarc[arrowsize=2pt
3,arrowinset=0]{->}(.6,.15){.15}{90}{210}\psarc(.6,.15){.15}{205}{270}
\end{pspicture} }

\newcommand{\SPCL}{\begin{pspicture}(0,.2)(.6,.7)
\psline(0,0)(.6,.3) \psline[arrowsize=2pt
3,arrowinset=0]{<-}(.1,.05)(.6,.3)\psarc[arrowsize=2pt
3,arrowinset=0]{->}(0,.45){.15}{270}{30}\psarc(0,.45){.15}{25}{90}
\psarc[arrowsize=2pt
3,arrowinset=0]{->}(.6,.3){.3}{90}{195}\psarc(.6,.3){.3}{195}{270}
\end{pspicture} }

\newcommand{\SPCM}{\begin{pspicture}(0,.2)(.6,.7)
\psline(0,.6)(.6,.3) \psline[arrowsize=2pt
3,arrowinset=0]{<-}(.1,.55)(.6,.3) \psarc[arrowsize=2pt
3,arrowinset=0]{->}(0,.15){.15}{270}{30}\psarc(0,.15){.15}{25}{90}
\psarc[arrowsize=2pt 3,
arrowinset=0]{->}(.6,.3){.3}{90}{210}\psarc(.6,.3){.3}{205}{270}
\end{pspicture} }

\newcommand{\SPCN}{\begin{pspicture}(0,.2)(.6,.7)
\psline(0,.3)(.6,.6) \psline[arrowsize=2pt 3,
arrowinset=0]{<-}(.05,.325)(.6,.6)\psarc[arrowsize=2pt 3,
arrowinset=0]{->}(0,.3){.3}{270}{15}\psarc(0,.3){.3}{15}{90}
\psarc[arrowsize=2pt 3,
arrowinset=0]{->}(.6,.15){.15}{90}{210}\psarc(.6,.15){.15}{205}{270}
\end{pspicture} }

\newcommand{\SPCO}{\begin{pspicture}(0,.2)(.6,.7)
\psline(0,.3)(.6,0) \psline[arrowsize=2pt 3,
arrowinset=0]{<-}(.05,.275)(.6,0) \psarc[arrowsize=2pt 3,
arrowinset=0]{->}(0,.3){.3}{270}{20}\psarc(0,.3){.3}{20}{90}
\psarc[arrowsize=2pt 3,
arrowinset=0]{->}(.6,.45){.15}{90}{210}\psarc(.6,.45){.15}{205}{270}
\end{pspicture} }

\newcommand{\eBA}{\begin{pspicture}(0,.2)(.8,.6)
\psline[arrowsize=2pt 3, arrowinset=0.75]{<-}(0,0)(.2,.3)
\psline[arrowsize=2pt 3, arrowinset=0.75]{<-}(0,.6)(.2,.3)
\psline[arrowsize=2pt 3,
arrowinset=0.75]{->}(.2,.3)(.6,.3)\psline(.2,.3)(.6,.3)
\psline[arrowsize=2pt 3,
arrowinset=0.75]{<-}(.6,.3)(.8,0)\psline(.6,.3)(.8,0)
\psline[arrowsize=2pt 3,
arrowinset=0.75]{<-}(.6,.3)(.8,.6)\psline(.6,.3)(.8,.6)
\end{pspicture} }

\newcommand{\eBAspec}{\begin{pspicture}(0,.2)(1,.6)
\psline[arrowsize=2pt 3, arrowinset=0.75]{<-}(0,0)(.2,.3)
\psline[arrowsize=2pt 3, arrowinset=0.75]{<-}(0,.6)(.2,.3)
\psline[arrowsize=2pt 3,
arrowinset=0.75]{->}(.2,.3)(.8,.3)\psline(.2,.3)(.8,.3)
\psline[arrowsize=2pt 3,
arrowinset=0.75]{<-}(.8,.3)(1,0)\psline(.8,.3)(1,0)
\psline[arrowsize=2pt 3,
arrowinset=0.75]{<-}(.8,.3)(1,.6)\psline(.8,.3)(1,.6)
\psline[arrowsize=2pt 6,arrowinset=.4]{<-}(.32,.3)(.6,.3)
\psline[dotsize=2pt 5.2]{o-}(.53,.3)(.6,.3) \psline(.4,.4)(.6,.4)
\end{pspicture} }

\newcommand{\eBB}{\begin{pspicture}(0,.2)(.8,.6)
\psline[arrowsize=2pt 3, arrowinset=0.75]{<-}(0,0)(.4,.1)
\psline[arrowsize=2pt 3,
arrowinset=0.75]{->}(0,.6)(.4,.5)\psline(0,.6)(.4,.5)
\psline[arrowsize=2pt 3,
arrowinset=0.75]{->}(.4,.1)(.4,.5)\psline(.4,.1)(.4,.5)
\psline[arrowsize=2pt 3,
arrowinset=0.75]{<-}(.4,.5)(.8,.6)\psline(.4,.5)(.8,.6)
\psline[arrowsize=2pt 3,
arrowinset=0.75]{->}(.4,.1)(.8,0)\psline(.4,.1)(.8,.0)
\end{pspicture} }

\newcommand{\eBC}{\begin{pspicture}(0,.2)(.8,.6)
\psline[arrowsize=2pt 3, arrowinset=0.75]{->}(0,0)(.4,.1)
\psline(0,0)(.4,.1)\psline[arrowsize=2pt 3,
arrowinset=0.75]{<-}(0,.6)(.4,.5)\psline(0,.6)(.4,.5)
\psline[arrowsize=2pt 3,
arrowinset=0.75]{<-}(.4,.1)(.4,.5)\psline(.4,.1)(.4,.5)
\psline[arrowsize=2pt 3,
arrowinset=0.75]{->}(.4,.5)(.8,0)\psline(.4,.5)(.8,0)
\psline[arrowsize=2pt 3,
arrowinset=0.75]{<-}(.4,.1)(.8,.6)\psline(.4,.1)(.8,.6)
\end{pspicture} }

\newcommand{\eBb}{\begin{pspicture}(0,.2)(.8,.6)
\psline[arrowsize=2pt 3, arrowinset=0.75]{->}(0,0)(.4,.1)
\psline(0,0)(.4,.1)\psline[arrowsize=2pt 3,
arrowinset=0.75]{<-}(0,.6)(.4,.5)\psline(0,.6)(.4,.5)
\psline[arrowsize=2pt 3,
arrowinset=0.75]{<-}(.4,.1)(.4,.5)\psline(.4,.1)(.4,.5)
\psline[arrowsize=2pt 3,
arrowinset=0.75]{->}(.4,.5)(.8,.6)\psline(.4,.5)(.8,.6)
\psline[arrowsize=2pt 3,
arrowinset=0.75]{<-}(.4,.1)(.8,0)\psline(.4,.1)(.8,.0)
\end{pspicture} }

\newcommand{\eCA}{\begin{pspicture}(0,.2)(.6,.6)
\psline[arrowsize=2pt 2.5, arrowinset=0.75]{<-}(0,.3)(.2,.45)
\psline[arrowsize=2pt 2.5, arrowinset=0.75]{<-}(0,.6)(.2,.45)
\psline(0,.3)(.2,.45) \psline(0,.6)(.2,.45)
\psline[arrowsize=2pt 2.5, arrowinset=0.75]{->}(.6,.3)(.4,.45)
\psline[arrowsize=2pt 2.5, arrowinset=0.75]{->}(.6,.6)(.4,.45)
\psline(.6,.3)(.4,.45) \psline(.6,.6)(.4,.45)

\psline[arrowsize=2pt 2.5,
arrowinset=0.75]{->}(.2,.45)(.4,.45)\psline(.2,.45)(.4,.45)
\psline[arrowsize=2pt 2.5, arrowinset=0.75]{<-}(0,0)(.6,0)
\psline(0,0)(.6,0)
\end{pspicture} }

\newcommand{\eCB}{\begin{pspicture}(0,.2)(.6,.6)
\psline[arrowsize=2pt 2.5, arrowinset=0.75]{<-}(0,.3)(.2,.45)
\psline[arrowsize=2pt 2.5, arrowinset=0.75]{<-}(0,.6)(.2,.45)
\psline(0,.3)(.2,.45) \psline(0,.6)(.2,.45)
\psline[arrowsize=2pt 2.5, arrowinset=0.75]{->}(.6,0)(.4,.3)
\psline[arrowsize=2pt 2.5, arrowinset=0.75]{->}(.6,.6)(.4,.3)
\psline(.6,0)(.4,.3) \psline(.6,.6)(.4,.3) \psline[arrowsize=2pt
2.5, arrowinset=0.75]{->}(.2,.45)(.4,.3)\psline(.2,.45)(.4,.3)
\psline[arrowsize=2pt 2.5,
arrowinset=0.75]{<-}(0,0)(.6,.3)\psline(.2,.45)(.4,.3)
\psline(0,0)(.6,.3)
\end{pspicture} }

\newcommand{\eCC}{\begin{pspicture}(0,.2)(.6,.6)
\psline[arrowsize=2pt 2.5, arrowinset=0.75]{<-}(0,.3)(.2,.45)
\psline[arrowsize=2pt 2.5, arrowinset=0.75]{<-}(0,.6)(.2,.45)
\psline(0,.3)(.2,.45) \psline(0,.6)(.2,.45)
\psline[arrowsize=2pt 2.5, arrowinset=0.75]{->}(.6,.3)(.4,.15)
\psline[arrowsize=2pt 2.5, arrowinset=0.75]{->}(.6,0)(.4,.15)
\psline(.6,.3)(.4,.15) \psline(.6,0)(.4,.15)

\psline[arrowsize=2pt 2.5,
arrowinset=0.75]{->}(.2,.45)(.4,.15)\psline(.2,.45)(.4,.15)
\psline[arrowsize=2pt 2.5, arrowinset=0.75]{<-}(0,0)(.6,.6)
\psline(0,0)(.6,.6)
\end{pspicture} }

\newcommand{\eCD}{\begin{pspicture}(0,.2)(.6,.6)
\psline[arrowsize=2pt 2.5, arrowinset=0.75]{<-}(0,0)(.2,.3)
\psline[arrowsize=2pt 2.5, arrowinset=0.75]{<-}(0,.6)(.2,.3)
\psline(0,0)(.2,.3) \psline(0,.6)(.2,.3)
\psline[arrowsize=2pt 2.5, arrowinset=0.75]{->}(.6,.3)(.4,.45)
\psline[arrowsize=2pt 2.5, arrowinset=0.75]{->}(.6,.6)(.4,.45)
\psline(.6,.3)(.4,.45) \psline(.6,.6)(.4,.45)

\psline[arrowsize=2pt 2.5,
arrowinset=0.75]{->}(.2,.3)(.4,.45)\psline(.2,.3)(.4,.45)
\psline[arrowsize=2pt 2.5, arrowinset=0.75]{<-}(0,.3)(.6,0)
\psline(0,.3)(.6,0)
\end{pspicture} }

\newcommand{\eCE}{\begin{pspicture}(0,.2)(.6,.6)
\psline[arrowsize=2pt 2.5, arrowinset=0.75]{<-}(0,0)(.2,.4)
\psline[arrowsize=2pt 2.5, arrowinset=0.75]{<-}(0,.6)(.2,.4)
\psline(0,0)(.2,.4) \psline(0,.6)(.2,.4)
\psline[arrowsize=2pt 2.5, arrowinset=0.75]{->}(.6,0)(.4,.4)
\psline[arrowsize=2pt 2.5, arrowinset=0.75]{->}(.6,.6)(.4,.4)
\psline(.6,0)(.4,.4) \psline(.6,.6)(.4,.4) \psline[arrowsize=2pt
2.5, arrowinset=0.75]{->}(.2,.4)(.4,.4)\psline(.2,.4)(.4,.4)
\psline[arrowsize=2pt 2.5, arrowinset=0.75]{<-}(0,.3)(.6,.3)
\psline(0,.3)(.6,.3)
\end{pspicture} }

\newcommand{\eCF}{\begin{pspicture}(0,.2)(.6,.6)
\psline[arrowsize=2pt 2.5, arrowinset=0.75]{<-}(0,0)(.2,.3)
\psline[arrowsize=2pt 2.5, arrowinset=0.75]{<-}(0,.6)(.2,.3)
\psline(0,0)(.2,.3) \psline(0,.6)(.2,.3)
\psline[arrowsize=2pt 2.5, arrowinset=0.75]{->}(.6,.3)(.4,.15)
\psline[arrowsize=2pt 2.5, arrowinset=0.75]{->}(.6,0)(.4,.15)
\psline(.6,.3)(.4,.15) \psline(.6,0)(.4,.15)

\psline[arrowsize=2pt 2.5,
arrowinset=0.75]{->}(.2,.3)(.4,.15)\psline(.2,.3)(.4,.15)
\psline[arrowsize=2pt 2.5, arrowinset=0.75]{<-}(0,.3)(.6,.6)
\psline(0,.3)(.6,.6)
\end{pspicture} }

\newcommand{\eCG}{\begin{pspicture}(0,.2)(.6,.6)
\psline[arrowsize=2pt 2.5, arrowinset=0.75]{<-}(0,.3)(.2,.15)
\psline[arrowsize=2pt 2.5, arrowinset=0.75]{<-}(0,0)(.2,.15)
\psline(0,.3)(.2,.15) \psline(0,0)(.2,.15)
\psline[arrowsize=2pt 2.5, arrowinset=0.75]{->}(.6,.3)(.4,.15)
\psline[arrowsize=2pt 2.5, arrowinset=0.75]{->}(.6,0)(.4,.15)
\psline(.6,.3)(.4,.15) \psline(.6,0)(.4,.15)

\psline[arrowsize=2pt 2.5,
arrowinset=0.75]{->}(.2,.15)(.4,.15)\psline(.2,.15)(.4,.15)
\psline[arrowsize=2pt 2.5, arrowinset=0.75]{<-}(0,.6)(.6,.6)
\psline(0,.6)(.6,.6)
\end{pspicture} }

\newcommand{\eCH}{\begin{pspicture}(0,.2)(.6,.6)
\psline[arrowsize=2pt 2.5, arrowinset=0.75]{<-}(0,.3)(.2,.15)
\psline[arrowsize=2pt 2.5, arrowinset=0.75]{<-}(0,0)(.2,.15)
\psline(0,.3)(.2,.15) \psline(0,0)(.2,.15)
\psline[arrowsize=2pt 2.5, arrowinset=0.75]{->}(.6,0)(.4,.3)
\psline[arrowsize=2pt 2.5, arrowinset=0.75]{->}(.6,.6)(.4,.3)
\psline(.6,0)(.4,.3) \psline(.6,.6)(.4,.3)

\psline[arrowsize=2pt 2.5,
arrowinset=0.75]{->}(.2,.15)(.4,.3)\psline(.2,.15)(.4,.3)
\psline[arrowsize=2pt 2.5, arrowinset=0.75]{<-}(0,.6)(.6,.3)
\psline(0,.6)(.6,.3)
\end{pspicture} }

\newcommand{\eCI}{\begin{pspicture}(0,.2)(.6,.6)
\psline[arrowsize=2pt 2.5, arrowinset=0.75]{<-}(0,.3)(.2,.15)
\psline[arrowsize=2pt 2.5, arrowinset=0.75]{<-}(0,0)(.2,.15)
\psline(0,.3)(.2,.15) \psline(0,0)(.2,.15)
\psline[arrowsize=2pt 2.5, arrowinset=0.75]{->}(.6,.3)(.4,.45)
\psline[arrowsize=2pt 2.5, arrowinset=0.75]{->}(.6,.6)(.4,.45)
\psline(.6,.3)(.4,.45) \psline(.6,.6)(.4,.45)

\psline[arrowsize=2pt 2.5,
arrowinset=0.75]{->}(.2,.15)(.4,.45)\psline(.2,.15)(.4,.45)
\psline[arrowsize=2pt 2.5, arrowinset=0.75]{<-}(0,.6)(.6,0)
\psline(0,.6)(.6,0)
\end{pspicture} }

\newcommand{\eCJ}{\begin{pspicture}(0,.2)(.6,.6)
\psline[arrowsize=2pt 2.5, arrowinset=0.75]{<-}(0,0)(.25,.3)
\psline[arrowsize=2pt 2.5,
arrowinset=0.75]{<-}(0,.6)(.25,.3)\psline[arrowsize=2pt 2.5,
arrowinset=0.75]{<-}(0,.3)(.25,.3) \psline(0,0)(.25,.3)
\psline(0,.6)(.25,.3)\psline(0,.3)(.25,.3)
\psline[arrowsize=2pt 2.5, arrowinset=0.75]{->}(.6,0)(.35,.3)
\psline[arrowsize=2pt 2.5,
arrowinset=0.75]{->}(.6,.6)(.35,.3)\psline[arrowsize=2pt 2.5,
arrowinset=0.75]{->}(.6,.3)(.4,.3) \psline(.6,0)(.35,.3)
\psline(.6,.6)(.35,.3)\psline(.6,.3)(.35,.3)
\end{pspicture} }

\newcommand{\eCK}{\begin{pspicture}(0,.2)(.6,.6)
\psline[arrowsize=2pt 2.5, arrowinset=0.75]{->}(.6,0)(.35,.3)
\psline[arrowsize=2pt 2.5,
arrowinset=0.75]{->}(.6,.6)(.35,.3)\psline[arrowsize=2pt 2.5,
arrowinset=0.75]{->}(.6,.3)(.4,.3) \psline(.6,0)(.35,.3)
\psline(.6,.6)(.35,.3)\psline(.6,.3)(.35,.3)
\end{pspicture} }

\newcommand{\eDA}{\begin{pspicture}(0,.35)(.9,.95)
\psline{<-}(0,.45)(.3,.45) \psline(.3,.45)(.6,.15)
\psline(.3,.45)(.6,.75)\psline(.6,.75)(.9,.9)
\psline(.6,.75)(.9,.6)\psline(.6,.15)(.9,.3) \psline(.6,.15)(.9,0)
\psline[arrowsize=2pt 2.5, arrowinset=0.75](0,.45)(.3,.45)
\psline[arrowsize=2pt 2.5, arrowinset=0.75]{->}(.3,.45)(.6,.15)
\psline[arrowsize=2pt 2.5,
arrowinset=0.75]{->}(.3,.45)(.6,.75)\psline[arrowsize=2pt 2.5,
arrowinset=0.75]{<-}(.6,.75)(.9,.9) \psline[arrowsize=2pt 2.5,
arrowinset=0.75]{<-}(.6,.75)(.9,.6)\psline[arrowsize=2pt 2.5,
arrowinset=0.75]{<-}(.6,.15)(.9,.3) \psline[arrowsize=2pt 2.5,
arrowinset=0.75]{<-}(.6,.15)(.9,0)
\end{pspicture} }

\newcommand{\eDB}{\begin{pspicture}(0,.35)(.7,.95)
\psline[arrowsize=2pt 2.5, arrowinset=0.75]{->}(.6,0)(.35,.3)
\psline[arrowsize=2pt 2.5,
arrowinset=0.75]{->}(.6,.6)(.35,.3)\psline[arrowsize=2pt 2.5,
arrowinset=0.75]{->}(.6,.3)(.4,.3) \psline(.6,0)(.35,.3)
\psline(.6,.6)(.35,.3)\psline(.6,.3)(.35,.3)
\psline(0,.4)(.6,.9)\psline[arrowsize=2pt 2.5,
arrowinset=0.75]{<-}(0,.4)(.6,.9)
\end{pspicture} }

\newcommand{\eeBA}{\begin{pspicture}(0,.2)(.6,.7)
\psline(0,0)(.6,.6) \psline(0,.6)(.6,0)\psline[arrowsize=2pt 4,
arrowinset=0]{<-}(.05,.05)(.6,.6) \psline[arrowsize=2pt 4,
arrowinset=0]{<-}(.05,.55)(.6,0)\pscircle*(.3,.3){.08}
\end{pspicture} }

\newcommand{\eeBAbs}{\begin{pspicture}(0,.1)(.4,.5)
\psline(0,0)(.4,.4) \psline(0,.4)(.4,0)\psline[arrowsize=2pt 4,
arrowinset=0]{<-}(.02,.02)(.4,.4) \psline[arrowsize=2pt 4,
arrowinset=0]{<-}(.23,.17)(.4,0)\pscircle*(.2,.2){.05}
\end{pspicture} }

\newcommand{\EEBA}{\begin{pspicture}(0,.2)(.6,.7)
\psline(0,0)(.6,.6) \psline(0,.6)(.6,0)\pscircle*(.3,.3){.08}
\end{pspicture} }

\newcommand{\eeCA}{\begin{pspicture}(0,.2)(.6,.7)
\pscircle*(.3,.45){.05} \psline(0,0)(.6,0) \psline(0,.3)(.6,.6)
\psline(0,.6)(.6,.3) \psline[arrowsize=2pt
3,arrowinset=0]{<-}(.05,.325)(.6,.6) \psline[arrowsize=2pt
3,arrowinset=0]{<-}(.05,.575)(.6,.3)
\end{pspicture} }

\newcommand{\eeCB}{\begin{pspicture}(0,.2)(.6,.7)
\pscircle*(.3,.15){.05} \psline(0,0)(.6,.3) \psline(0,.3)(.6,0)
\psline(0,.6)(.6,.6) \psline[arrowsize=2pt
3,arrowinset=0]{<-}(.05,.025)(.6,.3) \psline[arrowsize=2pt
3,arrowinset=0]{<-}(.05,.275)(.6,0)
\end{pspicture} }

\newcommand{\eeCC}{\begin{pspicture}(0,.2)(.6,.7)
\psarc(.3,.82){.6}{240}{300} \psarc(-.22,-.22){.85}{015}{075}
\psarc(.82,-.22){.85}{105}{165} \pscircle*(.3,.45){.05}
 \psarc[arrowsize=2pt
3,arrowinset=0]{->}(-.22,-.22){.85}{015}{070} \psarc[arrowsize=2pt
3,arrowinset=0]{->}(.82,-.22){.85}{105}{160}
\end{pspicture} }

\newcommand{\eeCD}{\begin{pspicture}(0,.2)(.6,.7)
\pscircle*(.4,.2){.05}\psline(0,0)(.6,.3) \psline(0,.3)(.6,.6)
\psline(0,.6)(.6,0) \psline[arrowsize=2pt
3,arrowinset=0]{<-}(.1,.05)(.6,.3)  \psline[arrowsize=2pt
3,arrowinset=0]{<-}(.225,.375)(.6,0)
\end{pspicture} }

\newcommand{\eeCE}{\begin{pspicture}(0,.2)(.6,.7)
\pscircle*(.2,.4){.05}\psline(0,0)(.6,.3) \psline(0,.3)(.6,.6)
\psline(0,.6)(.6,0)  \psline[arrowsize=2pt
3,arrowinset=0]{<-}(0,.3)(.6,.6) \psline[arrowsize=2pt
3,arrowinset=0]{<-}(0,.6)(.6,0)
\end{pspicture} }

\newcommand{\eeCF}{\begin{pspicture}(0,.2)(.6,.7)
\pscircle*(.4,.4){.05} \psline(0,0)(.6,.6) \psline(0,.3)(.6,0)
\psline(0,.6)(.6,.3) \psline[arrowsize=2pt
3,arrowinset=0]{<-}(.225,.225)(.6,.6)  \psline[arrowsize=2pt
3,arrowinset=0]{<-}(.1,.55)(.6,.3)
\end{pspicture} }

\newcommand{\eeCG}{\begin{pspicture}(0,.2)(.6,.7)
\pscircle*(.2,.2){.05} \psline(0,0)(.6,.6) \psline(0,.3)(.6,0)
\psline(0,.6)(.6,.3) \psline[arrowsize=2pt
3,arrowinset=0]{<-}(0,0)(.6,.6) \psline[arrowsize=2pt
3,arrowinset=0]{<-}(0,.3)(.6,0)
\end{pspicture} }

\newcommand{\eeCH}{\begin{pspicture}(0,.2)(.6,.7)
\psarc(.3,.82){.6}{240}{300} \psarc(-.22,-.22){.85}{015}{075}
\psarc(.82,-.22){.85}{105}{165} \pscircle*(.48,.25){.05}
 \psarc[arrowsize=2pt
3,arrowinset=0]{->}(-.22,-.22){.85}{015}{070} \psarc[arrowsize=2pt
3,arrowinset=0]{<-}(.3,.82){.6}{260}{300}
\end{pspicture} }

\newcommand{\eeCI}{\begin{pspicture}(0,.2)(.6,.7)
\psarc(.3,-.22){.6}{60}{120} \psarc(.82,.82){.85}{195}{255}
\psarc(-.22,.82){.85}{285}{345} \pscircle*(.48,.35){.05}
\psarc[arrowsize=2pt 3,arrowinset=0]{->}(.3,-.22){.6}{60}{100}
\psarc[arrowsize=2pt 3,arrowinset=0]{<-}(-.22,.82){.85}{290}{345}
\end{pspicture} }

\newcommand{\eeCJ}{\begin{pspicture}(0,.2)(.6,.7)
\pscircle*(.3,.3){.05}  \psarc[arrowsize=2pt
3,arrowinset=0]{->}(0,.45){.15}{270}{30}\psarc(0,.45){.15}{30}{90}
\psline(0,0)(.6,.6)\psline[arrowsize=2pt
3,arrowinset=0]{<-}(.1,.1)(.6,.6)\psline(.3,.3)(.6,0)\psline(.3,.3)(.6,.3)
\end{pspicture} }

\newcommand{\eeCK}{\begin{pspicture}(0,.2)(.6,.7)
\pscircle*(.3,.3){.05} \psarc[arrowsize=2pt
3,arrowinset=0]{->}(0,.15){.15}{270}{30}\psarc(0,.15){.15}{25}{90}
\psline(0,.6)(.6,0)\psline[arrowsize=2pt
3,arrowinset=0]{<-}(.1,.5)(.6,0)\psline(.3,.3)(.6,.6)\psline(.3,.3)(.6,.3)
\end{pspicture} }

\newcommand{\eeCL}{\begin{pspicture}(0,.2)(.6,.7)
\pscircle*(.3,.3){.05}
\psline(0,.3)(.3,.3)\psline(0,0)(.3,.3)\psline(0,.6)(.6,0)
\psline[arrowsize=2pt
3,arrowinset=0]{<-}(0,.3)(.3,.3)\psline[arrowsize=2pt
3,arrowinset=0]{<-}(.1,.1)(.3,.3)\psline[arrowsize=2pt
3,arrowinset=0]{<-}(.1,.5)(.3,.3)\psarc[arrowsize=2pt
3,arrowinset=0]{->}(.6,.45){.15}{90}{210}\psarc(.6,.45){.15}{210}{270}
\end{pspicture} }

\newcommand{\eeCM}{\begin{pspicture}(0,.2)(.6,.7)
\pscircle*(.3,.3){.05}\psline(0,.6)(.3,.3)
\psline(0,.3)(.3,.3)\psline(0,0)(.6,.6) \psline[arrowsize=2pt
3,arrowinset=0]{<-}(.1,.5)(.3,.3) \psline[arrowsize=2pt
3,arrowinset=0]{<-}(0,.3)(.3,.3)\psline[arrowsize=2pt
3,arrowinset=0]{<-}(.1,.1)(.6,.6)\psarc[arrowsize=2pt
3,arrowinset=0]{->}(.6,.15){.15}{90}{210}\psarc(.6,.15){.15}{205}{270}
\end{pspicture} }

\newcommand{\eeCN}{\begin{pspicture}(0,.2)(.7,.7)
\pscircle*(.5,.3){.05} \psline(0,.3)(.7,.3)\psline(0,.3)(.7,.3)
\psline[arrowsize=2pt
3,arrowinset=0]{<-}(.05,.3)(.7,.3)\psline(.5,.3)(.7,.6)
\psline(.5,.3)(.7,0) \psarc[arrowsize=2pt
3,arrowinset=0]{->}(0,.3){.3}{270}{70}\psarc(0,.3){.3}{270}{90}
\end{pspicture} }

\newcommand{\eeCO}{\begin{pspicture}(0,.2)(.7,.7)
\pscircle*(.2,.3){.05} \psline(0,.3)(.7,.3) \psline(.2,.3)(0,.6)
\psline(.2,.3)(0,0)\psline[arrowsize=2pt
3,arrowinset=0]{<-}(-.1,.3)(.7,.3) \psline[arrowsize=2pt
3,arrowinset=0]{->}(.2,.3)(0,.6) \psline[arrowsize=2pt
3,arrowinset=0]{->}(.2,.3)(0,0) \psarc[arrowsize=2pt
3,arrowinset=0]{->}(.7,.3){.3}{90}{250} \psarc(.7,.3){.3}{90}{270}
\end{pspicture} }

\newcommand{\eeCP}{\begin{pspicture}(0,.2)(.6,.7)
\psarc(.3,.80){.6}{230}{310} \psarc(-.22,-.22){.85}{015}{075}
\psarc(.82,-.22){.85}{105}{165} \psarc[arrowsize=2pt
3,arrowinset=0]{<-}(.3,.80){.6}{230}{310} \psarc[arrowsize=2pt
3,arrowinset=0]{->}(-.22,-.22){.85}{015}{075} \psarc[arrowsize=2pt
3,arrowinset=0]{->}(.82,-.22){.85}{105}{165} \psarc[arrowsize=2pt
3,arrowinset=0]{->}(.82,-.22){.85}{105}{144}
\pscircle*(.3,.44){.05}\pscircle*(.11,.24){.05}
\end{pspicture} }

\newcommand{\eeCQ}{\begin{pspicture}(0,.2)(.6,.7)
\psarc(.3,.80){.6}{230}{310} \psarc(-.22,-.22){.85}{015}{075}
\psarc(.82,-.22){.85}{105}{165}\psarc[arrowsize=2pt
3,arrowinset=0]{<-}(.3,.80){.6}{230}{310} \psarc[arrowsize=2pt
3,arrowinset=0]{->}(-.22,-.22){.85}{015}{075}\psarc[arrowsize=2pt
3,arrowinset=0]{->}(-.22,-.22){.85}{015}{048} \psarc[arrowsize=2pt
3,arrowinset=0]{->}(.82,-.22){.85}{105}{165}
\pscircle*(.3,.44){.05} \pscircle*(.49,.24){.05}
\end{pspicture} }

\newcommand{\eeCR}{\begin{pspicture}(0,.2)(.6,.7)
\psarc(.3,.80){.6}{230}{310} \psarc(-.22,-.22){.85}{015}{075}
\psarc(.82,-.22){.85}{105}{165}\psarc[arrowsize=2pt
3,arrowinset=0]{<-}(.3,.80){.6}{230}{310} \psarc[arrowsize=2pt
3,arrowinset=0]{<-}(.3,.80){.6}{260}{310} \psarc[arrowsize=2pt
3,arrowinset=0]{->}(-.22,-.22){.85}{015}{075} \psarc[arrowsize=2pt
3,arrowinset=0]{->}(.82,-.22){.85}{105}{165}\pscircle*(.11,.24){.05}
\pscircle*(.49,.24){.05}
\end{pspicture} }

\newcommand{\eeCS}{\begin{pspicture}(0,.2)(.6,.7)
\psarc(.3,-.2){.6}{50}{130} \psarc(.82,.82){.85}{195}{255}
\psarc(-.22,.82){.85}{285}{345}\psarc[arrowsize=2pt
3,arrowinset=0]{->}(.3,-.2){.6}{50}{130} \psarc[arrowsize=2pt
3,arrowinset=0]{<-}(.82,.82){.85}{195}{255}\psarc[arrowsize=2pt
3,arrowinset=0]{<-}(.82,.82){.85}{215}{255}  \psarc[arrowsize=2pt
3,arrowinset=0]{<-}(-.22,.82){.85}{285}{345}\pscircle*(.3,.15){.05}
\pscircle*(.11,.37){.05}
\end{pspicture} }

\newcommand{\eeCT}{\begin{pspicture}(0,.2)(.6,.7)
\psarc(.3,-.2){.6}{50}{130} \psarc(.82,.82){.85}{195}{255}
\psarc(-.22,.82){.85}{285}{345}\psarc[arrowsize=2pt
3,arrowinset=0]{->}(.3,-.2){.6}{50}{130} \psarc[arrowsize=2pt
3,arrowinset=0]{<-}(.82,.82){.85}{195}{255} \psarc[arrowsize=2pt
3,arrowinset=0]{<-}(-.22,.82){.85}{285}{345}\psarc[arrowsize=2pt
3,arrowinset=0]{<-}(-.22,.82){.85}{312}{345}\pscircle*(.3,.15){.05}
\pscircle*(.49,.37){.05}
\end{pspicture} }

\newcommand{\eeCU}{\begin{pspicture}(0,.2)(.6,.7)
\psarc(.3,-.2){.6}{50}{130} \psarc(.82,.82){.85}{195}{255}
\psarc(-.22,.82){.85}{285}{345}\psarc[arrowsize=2pt
3,arrowinset=0]{->}(.3,-.2){.6}{50}{130}\psarc[arrowsize=2pt
3,arrowinset=0]{->}(.3,-.2){.6}{50}{100} \psarc[arrowsize=2pt
3,arrowinset=0]{<-}(.82,.82){.85}{195}{255} \psarc[arrowsize=2pt
3,arrowinset=0]{<-}(-.22,.82){.85}{285}{345}\pscircle*(.11,.37){.05}
\pscircle*(.49,.37){.05}
\end{pspicture} }

\newcommand{\eeCV}{\begin{pspicture}(0,.2)(.6,.7)
\pscircle*(.2,.4){.05}\psline(0,0)(.6,.3) \psline(0,.3)(.6,.6)
\psline(0,.6)(.6,0)  \psline[arrowsize=2pt
3,arrowinset=0]{<-}(0,.3)(.6,.6) \psline[arrowsize=2pt
3,arrowinset=0]{<-}(0,.6)(.6,0) \pscircle*(.4,.2){.05}
\psline[arrowsize=2pt 3,arrowinset=0]{<-}(.1,.05)(.6,.3)
\psline[arrowsize=2pt 3,arrowinset=0]{<-}(.225,.375)(.6,0)
\end{pspicture} }

\newcommand{\eeCW}{\begin{pspicture}(0,.2)(.6,.7)
\pscircle*(.4,.4){.05} \psline(0,0)(.6,.6) \psline(0,.3)(.6,0)
\psline(0,.6)(.6,.3) \psline[arrowsize=2pt
3,arrowinset=0]{<-}(.225,.225)(.6,.6)  \psline[arrowsize=2pt
3,arrowinset=0]{<-}(.1,.55)(.6,.3)\pscircle*(.2,.2){.05}
\psline[arrowsize=2pt 3,arrowinset=0]{<-}(0,0)(.6,.6)
\psline[arrowsize=2pt 3,arrowinset=0]{<-}(0,.3)(.6,0)
\end{pspicture} }

\newcommand{\eeCX}{\begin{pspicture}(0,.2)(.7,.7)
\pscircle*(.2,.3){.05}\pscircle*(.5,.3){.05} \psline(0,.3)(.8,.3)
\psline(.2,.3)(0,.6) \psline(.2,.3)(0,0)\psline[arrowsize=2pt
3,arrowinset=0]{<-}(-.1,.3)(.7,.3)\psline[arrowsize=2pt
3,arrowinset=0]{<-}(.25,.3)(.7,.3) \psline[arrowsize=2pt
3,arrowinset=0]{->}(.2,.3)(0,.6) \psline[arrowsize=2pt
3,arrowinset=0]{->}(.2,.3)(0,0)\psline(.5,.3)(.7,.6)
\psline(.5,.3)(.7,0)
\end{pspicture} }

\newcommand{\eeCY}{\begin{pspicture}(0,.2)(.6,.7)
\pscircle*(.4,.4){.05} \psline(0,0)(.6,.6) \psline[arrowsize=2pt
3,arrowinset=0]{<-}(.225,.225)(.6,.6)\pscircle*(.2,.2){.05}
\psline[arrowsize=2pt 3,arrowinset=0]{<-}(0,0)(.6,.6)
\psline(0,.6)(.2,.2) \psline[arrowsize=2pt
3,arrowinset=0]{<-}(0,.6)(.2,.2)\psline(.4,.4)(.6,0)
\psline(-.1,.3)(.4,.4)\psline[arrowsize=2pt
3,arrowinset=0]{<-}(-.1,.3)(.4,.4) \psline(.2,.2)(.7,.3)
\end{pspicture} }

\newcommand{\eeDA}{\begin{pspicture}(0,.2)(.7,.7)
\pscircle*(.4,.3){.05} \psline(.4,.3)(.7,0) \psline(.4,.3)(.7,.2)
\psline(.4,.3)(.7,.4) \psline(.4,.3)(.7,.6)
\end{pspicture} }

\newcommand{\EECA}{\begin{pspicture}(0,.2)(.6,.7)
\pscircle*(.3,.45){.05} \psline(0,0)(.6,0) \psline[arrowsize=2pt
3,arrowinset=0]{<-}(.2,0)(.6,0)\psline(0,.3)(.6,.6)
\psline(0,.6)(.6,.3)
\end{pspicture} }

\newcommand{\EECB}{\begin{pspicture}(0,.2)(.6,.7)
\pscircle*(.3,.15){.05} \psline(0,0)(.6,.3)\psline(0,0)(.6,.3)
\psline(0,.3)(.6,0) \psline[arrowsize=2pt
3,arrowinset=0]{<-}(.2,.6)(.6,.6)\psline(0,.6)(.6,.6)
\end{pspicture} }

\newcommand{\EECD}{\begin{pspicture}(0,.2)(.6,.7)
\pscircle*(.4,.2){.05}\psline(0,0)(.6,.3)
\psline(0,.3)(.6,.6)\psline[arrowsize=2pt
3,arrowinset=0]{<-}(0,.3)(.6,.6) \psline(0,.6)(.6,0)
\end{pspicture} }

\newcommand{\EECE}{\begin{pspicture}(0,.2)(.6,.7)
\pscircle*(.2,.4){.05}\psline(0,0)(.6,.3)\psline[arrowsize=2pt
3,arrowinset=0]{<-}(.1,.05)(.6,.3) \psline(0,.3)(.6,.6)
\psline(0,.6)(.6,0)
\end{pspicture} }

\newcommand{\EECF}{\begin{pspicture}(0,.2)(.6,.7)
\pscircle*(.4,.4){.05} \psline(0,0)(.6,.6)
\psline(0,.3)(.6,0)\psline[arrowsize=2pt
3,arrowinset=0]{<-}(0,.3)(.6,0) \psline(0,.6)(.6,.3)
\end{pspicture} }

\newcommand{\EECG}{\begin{pspicture}(0,.2)(.6,.7)
\pscircle*(.2,.2){.05} \psline(0,0)(.6,.6) \psline(0,.3)(.6,0)
\psline(0,.6)(.6,.3)\psline[arrowsize=2pt
3,arrowinset=0]{<-}(.1,.55)(.6,.3)
\end{pspicture} }

\newcommand{\EECH}{\begin{pspicture}(0,.2)(.6,.7)
\psarc(.3,.82){.6}{240}{300} \psarc(-.22,-.22){.85}{015}{075}
\psarc(.82,-.22){.85}{105}{165} \psarc[arrowsize=2pt
3,arrowinset=0]{->}(.82,-.22){.85}{105}{145}\pscircle*(.48,.25){.05}
\end{pspicture} }

\newcommand{\EECI}{\begin{pspicture}(0,.2)(.6,.7)
\psarc(.3,-.22){.6}{60}{120} \psarc(.82,.82){.85}{195}{255}
\psarc[arrowsize=2pt 3,arrowinset=0]{<-}(.82,.82){.85}{215}{235}
\psarc(-.22,.82){.85}{285}{345} \pscircle*(.48,.35){.05}

\end{pspicture} }

\newcommand{\EECJ}{\begin{pspicture}(0,.2)(.6,.7)
\pscircle*(.3,.3){.05}  \psarc[arrowsize=2pt
3,arrowinset=0]{->}(0,.45){.15}{270}{30}\psarc(0,.45){.15}{30}{90}
\psline(0,0)(.6,.6)\psline(.3,.3)(.6,0)\psline(.3,.3)(.6,.3)
\end{pspicture} }

\newcommand{\EECK}{\begin{pspicture}(0,.2)(.6,.7)
\pscircle*(.3,.3){.05} \psarc[arrowsize=2pt
3,arrowinset=0]{->}(0,.15){.15}{270}{30}\psarc(0,.15){.15}{25}{90}
\psline(0,.6)(.6,0)\psline(.3,.3)(.6,.6)\psline(.3,.3)(.6,.3)
\end{pspicture} }

\newcommand{\EECL}{\begin{pspicture}(0,.2)(.6,.7)
\pscircle*(.3,.3){.05}
\psline(0,.3)(.3,.3)\psline(0,0)(.3,.3)\psline(0,.6)(.6,0)\psarc[arrowsize=2pt
3,arrowinset=0]{->}(.6,.45){.15}{90}{210}\psarc(.6,.45){.15}{210}{270}
\end{pspicture} }

\newcommand{\EECM}{\begin{pspicture}(0,.2)(.6,.7)
\pscircle*(.3,.3){.05}\psline(0,.6)(.3,.3)
\psline(0,.3)(.3,.3)\psline(0,0)(.6,.6)\psarc[arrowsize=2pt
3,arrowinset=0]{->}(.6,.15){.15}{90}{210}\psarc(.6,.15){.15}{205}{270}
\end{pspicture} }

\newcommand{\EECN}{\begin{pspicture}(0,.2)(.7,.7)
\pscircle*(.5,.3){.05} \psline(0,.3)(.7,.3) \psline(.5,.3)(.7,.6)
\psline(.5,.3)(.7,0) \psarc[arrowsize=2pt
3,arrowinset=0]{->}(0,.3){.3}{270}{70}\psarc(0,.3){.3}{270}{90}
\end{pspicture} }

\newcommand{\EECO}{\begin{pspicture}(0,.2)(.7,.7)
\pscircle*(.2,.3){.05} \psline(0,.3)(.7,.3) \psline(.2,.3)(0,.6)
\psline(.2,.3)(0,0) \psarc[arrowsize=2pt
3,arrowinset=0]{->}(.7,.3){.3}{90}{250} \psarc(.7,.3){.3}{90}{270}
\end{pspicture} }

\newcommand{\EECP}{\begin{pspicture}(0,.2)(.6,.7)
\psarc(.3,.80){.6}{230}{310} \psarc(-.22,-.22){.85}{015}{075}
\psarc(.82,-.22){.85}{105}{165}  \psarc[arrowsize=2pt
3,arrowinset=0]{->}(.82,-.22){.85}{105}{144}
\pscircle*(.3,.44){.05}\pscircle*(.11,.24){.05}
\end{pspicture} }

\newcommand{\EECQ}{\begin{pspicture}(0,.2)(.6,.7)
\psarc(.3,.80){.6}{230}{310} \psarc(-.22,-.22){.85}{015}{075}
\psarc(.82,-.22){.85}{105}{165}\psarc[arrowsize=2pt
3,arrowinset=0]{->}(-.22,-.22){.85}{015}{048}
\pscircle*(.3,.44){.05} \pscircle*(.49,.24){.05}
\end{pspicture} }

\newcommand{\EECR}{\begin{pspicture}(0,.2)(.6,.7)
\psarc(.3,.80){.6}{230}{310} \psarc(-.22,-.22){.85}{015}{075}
\psarc(.82,-.22){.85}{105}{165}\psarc[arrowsize=2pt
3,arrowinset=0]{<-}(.3,.80){.6}{260}{310}\pscircle*(.11,.24){.05}
\pscircle*(.49,.24){.05}
\end{pspicture} }

\newcommand{\EECS}{\begin{pspicture}(0,.2)(.6,.7)
\psarc(.3,-.2){.6}{50}{130} \psarc(.82,.82){.85}{195}{255}
\psarc(-.22,.82){.85}{285}{345}\psarc[arrowsize=2pt
3,arrowinset=0]{<-}(.82,.82){.85}{215}{255}
\pscircle*(.3,.15){.05} \pscircle*(.11,.37){.05}
\end{pspicture} }

\newcommand{\EECT}{\begin{pspicture}(0,.2)(.6,.7)
\psarc(.3,-.2){.6}{50}{130} \psarc(.82,.82){.85}{195}{255}
\psarc(-.22,.82){.85}{285}{345}\psarc[arrowsize=2pt
3,arrowinset=0]{<-}(-.22,.82){.85}{312}{345}\pscircle*(.3,.15){.05}
\pscircle*(.49,.37){.05}
\end{pspicture} }

\newcommand{\EECU}{\begin{pspicture}(0,.2)(.6,.7)
\psarc(.3,-.2){.6}{50}{130} \psarc(.82,.82){.85}{195}{255}
\psarc(-.22,.82){.85}{285}{345}\psarc[arrowsize=2pt
3,arrowinset=0]{->}(.3,-.2){.6}{50}{100} \pscircle*(.11,.37){.05}
\pscircle*(.49,.37){.05}
\end{pspicture} }

\newcommand{\EECV}{\begin{pspicture}(0,.2)(.6,.7)
\pscircle*(.2,.4){.05}\psline(0,0)(.6,.3) \psline(0,.3)(.6,.6)
\psline(0,.6)(.6,0)    \pscircle*(.4,.2){.05}
\psline[arrowsize=2pt 3,arrowinset=0]{<-}(.225,.375)(.6,0)
\end{pspicture} }

\newcommand{\EECW}{\begin{pspicture}(0,.2)(.6,.7)
\pscircle*(.4,.4){.05} \psline(0,0)(.6,.6) \psline(0,.3)(.6,0)
\psline(0,.6)(.6,.3) \psline[arrowsize=2pt
3,arrowinset=0]{<-}(.225,.225)(.6,.6) \pscircle*(.2,.2){.05}
\end{pspicture} }

\newcommand{\EECZ}{\begin{pspicture}(0,.2)(.6,.7)
\psarc(.3,.80){.6}{230}{310} \psarc(-.22,-.22){.85}{015}{075}
\psarc(.82,-.22){.85}{105}{165}  \psarc[arrowsize=2pt
3,arrowinset=0]{->}(.82,-.22){.85}{105}{144}
\pscircle*(.3,.44){.05}\pscircle*(.11,.24){.05}
 \pscircle*(.49,.24){.05}\psarc[arrowsize=2pt
3,arrowinset=0]{->}(-.22,-.22){.85}{015}{048} \psarc[arrowsize=2pt
3,arrowinset=0]{<-}(.3,.80){.6}{260}{310}
\end{pspicture} }

\newcommand{\EECZz}{\begin{pspicture}(0,.2)(.6,.7)
\psarc(.3,-.2){.6}{50}{130} \psarc(.82,.82){.85}{195}{255}
\psarc(-.22,.82){.85}{285}{345}\psarc[arrowsize=2pt
3,arrowinset=0]{<-}(.82,.82){.85}{215}{255}
\pscircle*(.3,.15){.05} \pscircle*(.11,.37){.05}
\pscircle*(.49,.37){.05} \psarc[arrowsize=2pt
3,arrowinset=0]{<-}(-.22,.82){.85}{312}{345} \psarc[arrowsize=2pt
3,arrowinset=0]{->}(.3,-.2){.6}{50}{100}
\end{pspicture} }
\makeatletter

\@addtoreset{equation}{section}

\def\section{\@startsection {section}{1}{\z@}{-3.5ex plus -1ex minus
 -.2ex}{2.3ex plus .2ex}{\large\bf\centering}}
\def\subsection{\@startsection{subsection}{2}{\z@}{-3.25ex plus%
 -1ex minus -.2ex}{1.5ex plus .2ex}{\bf}}
\def\subsubsection{\@startsection{subsubsection}{3}{\z@}{-3.25ex plus%
 -1ex minus -.2ex}{1.5ex plus .2ex}{\sl}}
\makeatother
\begin{document}

\baselineskip 19pt

\parindent 10pt
\parskip 17pt

\begin{titlepage}

\vspace{2.5cm}
\begin{center}
{\large {\bf Rational $R$-matrices, centralizer algebras and
tensor
identities\\[0.1in] for $e_6$ and $e_7$ exceptional families of Lie algebras}}\\
\vspace{1cm} {\large N. J. MacKay\footnote{\tt nm15@york.ac.uk}
and A. Taylor\footnote{\tt at165@york.ac.uk}}
\\
\vspace{3mm} {\em Department of Mathematics,\\ University of York,
\\York YO10 5DD, U.K.}
\end{center}

\vspace{0.5cm}
\begin{abstract}
\noindent We use Cvitanovi\'c's diagrammatic techniques to
construct the rational solutions of the Yang-Baxter equation
associated with the $e_6$ and $e_7$ families of Lie algebras, and
thus explain Westbury's observations about their uniform spectral
decompositions. In doing so we explore the extensions of the
Brauer and symmetric group algebras to the centralizer algebras of
$e_7$ and $e_6$ on their lowest-dimensional representations and
(up to three-fold) tensor products thereof, giving bases for them
and a range of identities satisfied by the algebras' defining
invariant tensors.
\end{abstract}

\end{titlepage}

\section{Introduction}

The Yang-Baxter equation (YBE) \cite{YBE}, which appears in $1+1D$
physics as the factorizability condition for $S$-matrices in
integrable models, is closely bound up with Lie algebras and their
representation theory, essentially because of the asymptotic
behaviour of its rational solutions (see (\ref{basicr}) below).
Indeed, if one were to investigate the YBE knowing nothing of Lie
algebras, one would very soon find oneself re-discovering a great
deal about them. In fact surprisingly little is known about the
rational YBE solutions associated with the exceptional Lie
algebras: this paper investigates these, and finds, in precisely
the spirit of the preceding sentence, an intricate relationship
between the YBE and various identities satisfied by the algebras'
invariant tensors.

Our point of departure is a remarkable observation made a few
years ago by Westbury \cite{westb02}: that certain solutions of
the YBE (`$R$-matrices') associated with the Lie algebras of the
$e_6$ and $e_7$ series (the second and third rows of the
Freudenthal-Tits `magic square') have spectral decompositions
which may be expressed simply and uniformly in terms of the
dimension ($=1,2,4$ or $8$) of the underlying division algebra. In
this paper we shall explicitly construct these $R$-matrices, prove
that they solve the Yang-Baxter equation, and thus provide an
explanation of Westbury's observation. More interesting, perhaps,
is what we shall learn along the way. In particular we will need
to understand the structure of, and provide a basis for, the
centralizer of the Lie group action on tensor cubes of the
defining representation. (These are the analogues for the
exceptional series of the symmetric group algebra for $su(n)$ and
of Brauer's algebra for the other classical groups.) We shall also
discover a host of secondary identities satisfied by the groups'
defining invariant tensors, all of them subtly necessary in
solving the YBE.

Our method
 is to use Cvitanovi\'c's `birdtrack' diagrams \cite{cvitabk,cvita76}
(which extend earlier ideas of Penrose) to handle the
calculations. An alternative approach to the $e_n$ centralizers,
which utilizes the braid matrices (the $q$-deformed but
spectral-parameter $u$-independent $R$-matrices) and is
complementary to ours, appears in \cite{wenzl03}.

The paper is structured as follows. In section two we provide a
brief recapitulation of some of Westbury's observations. In
section three we give an elementary, essentially pedagogical
recapitulation of these issues for the classical groups -- the
rational $R$-matrices, the centralizer algebras and the
diagrammatic techniques used to handle them. Section four deals
with the $e_6$ series, and section five with the $e_7$ series.

\section{The Yang-Baxter equation and a unified spectral decomposition
for exceptional R-matrices}

\subsection{The Yang-Baxter equation}

The Yang-Baxter equation (YBE), between expressions in
$End(V\otimes V \otimes V)$ for $V=\C^n$, is \be\label{YBE1}
\check{R}(u)\otimes 1\cdot 1\otimes \check{R}(u+v)\cdot
\check{R}(v)\otimes 1 = 1\otimes \check{R}(v)\cdot
\check{R}(u+v)\otimes 1\cdot1\otimes \check{R}(u)\,,\ee or, with
its indices made explicit (each running from $1$ to $n$, and with
repeated indices summed),\be\label{YBE2} \check{R}^{ij}_{lm}(u)
\check{R}^{mk}_{sr}(u+v) \check{R}^{ls}_{pq}(v) =
\check{R}_{lm}^{jk}(v) \check{R}_{ps}^{il}(u+v)
\check{R}_{qr}^{sm}(u)\,,\ee for $\check{R}(u)\in End(V\otimes
V)$. We first note that this equation is homogeneous in
$\check{R}$ and in $u$, so that $\mu \check{R}(\lambda  u)$ is
still a solution for arbitrary $\C$-scalings $\lambda$ and $\mu$.
We shall therefore rescale both $\check{R}$ and $u$ wherever it is
convenient for us to do so. (In the physical construction of
factorized $S$-matrices, in contrast, the scale of $u$ is fixed,
and scaling $\check{R}$ affects its analytic properties and thus
the bootstrap spectrum.)

The simplest class of solutions of the YBE (which we refer to as
`$R$-matrices') has rational dependence on $u$, and an expansion
in powers of $1/u$ of the form \be\label{basicr} \check{R}(u)
=\transp\left( \I_n\otimes \I_n + {C\over u}+\ldots\right)
\hspace{0.3in}{\rm
where}\hspace{0.3in}C=\sum_{a,b}\rho_V(I^a)\otimes\rho_V(I^b)g_{ab}\,,\ee
in which $\I_n$ is the $n\times n$ identity matrix, $I^a$ are the
generators of a Lie algebra $\g$, $g_{ab}$ its Cartan-Killing
form, $\rho_V$ its suitably-chosen representation on a module $V$
(usually its defining representation), and $\transp$ the
transposition operator on the two components of $V\otimes V$.
Thus, from the outset, the investigation of $R$-matrices naturally
involves the investigation of Lie algebras and their
representations.

A natural consequence of this (see, for example, \cite{macka91})
is that $\check{R}(u)$ commutes with the action of $\g$ on
$V\otimes V$, so that, by Schur's lemma, \be\label{spectral}
\check{R}(u) = \sum_i f_i(u)P_i \ee for some scalar functions
$f_i(u)$, where the sum is over projectors onto irreducible
components $W_i\subset V\otimes V$. (This is only fully correct
where there are no multiplicities. Where such repetitions among
the $W_i$ occur, there can be non-trivial intertwiners between
them.)

\subsection{$R$-matrix spectra and the magic square}

Now recall that the Freudenthal-Tits `magic square'
\cite{tits62,freud65} is
$$\begin{array}{c|cccc} m= & \quad 1 & 2 & 4\quad & 8 \\ \hline
&\quad  a_1 & a_2 & c_3\quad  & f_4
\\ & \quad a_2 & a_2 \times a_2 & a_5 \quad & e_6 \\
 & \quad c_3 & a_5 & d_6\quad  & e_7
\\ &\quad  f_4 & e_6 & e_7\quad  & e_8
\end{array}$$
(We will not need the details of its construction. For full
discussions, including an explanation of its
row$\leftrightarrow$column symmetry, see \cite{barto02,ramond}.)
We will refer to the row whose last ($m=8$) entry is the
exceptional algebra $\g$ as the `$\g$ series' of Lie algebras.

For the $e_6$ series, Westbury's principal observation in
\cite{westb02} was that, in the literature of rational $R$-matrix
spectral decompositions for individual $\g$ and $V$ (originally in
\cite{kulis81} for $a_n$, \cite{ogiev87} for $e_6$), there is a
unified underlying formula: for the representation on $V$ of
dimension $n=3m+3$, \be\label{e6R} \check{R}(u) =P_1 + {4+u\over
4-u} P_2 + {4+u\over 4-u}\; {2m+u\over 2m-u}P_3\,,\ee where $W_1$
is the representation whose highest weight is double that of $V$,
$W_2$ is the antisymmetric component of $V\otimes V$, and
$W_3=\bar{V}$, the complex-conjugate of $V$.

The YBE is straightforwardly generalized to act on $V_1\otimes
V_2\otimes V_3$ for $V_1\neq V_2 \neq V_3$. There is then a
unified spectral decomposition for
$\transp\check{R}_{V\bar{V}}(u)\in End(V\otimes \bar{V})$ (in
which $\transp$ now transposes elements of $V\otimes\bar{V}$ with
those of $\bar{V}\otimes V$), for which \be\label{e6Rcr}
\check{R}_{V\bar{V}}(u)= \transp\left(P_1 + {u+m+4\over u-m-4}P_2
+ {u+m+4\over u-m-4}\;{u+3m\over u-3m}P_3\right)\,,\ee where $W_1$
is the representation whose highest weight is the sum of those of
$V$ and $\bar{V}$, $W_2=\g$, the adjoint representation, and
$W_3=\C$, the singlet.

For the $e_7$ series, Westbury observes that, for $V$ of dimension
$n=6m+8$, \be\label{e7R} \check{R}(u) = P_1 + {2+u\over 2-u} P_2 +
{2+u\over 2-u} {m+2+u\over m+2-u} P_3 + {2+u\over 2-u}\;
{m+2+u\over m+2-u}{2m+2+u\over 2m+2-u} P_4\,,\ee where the highest
weight of $W_1$ is twice that of $V$, $W_2$ is the highest
antisymmetric component of $V\otimes V$, $W_3=\g$ and $W_4=\C$.
(The original $R$-matrix spectra are in \cite{macka92} for $c_3$,
\cite{kulis81} for $a_5$, \cite{shank79} for $d_6$ and
\cite{ogiev87} for $e_7$; see also \cite{westb03} for an extension
to further values of $m$.)

We shall not, in this paper, concern ourselves with the $g_2$
series (the `zeroth' row of the magic square, for which the
$R$-matrices are dealt with in \cite{ogiev86,macka04}), or the
$f_4$ and $e_8$ series, which are each, in different ways,
problematic.

For the $f_4$ series, where the same observation might be expected
to hold for $V$ of dimension $3m+2$, in fact (surprisingly) it
fails. A uniform decomposition exists for $c_3$ and $f_4$, but
fails to work fully for the other algebras in the series. We
suspect that the resolution is bound up with the identities
satisfied by the primitive invariant tensor, and are working to
understand this. A common feature of the $f_4$ and $e_8$
calculations is the need to evaluate `pentagon' diagrams (in the
diagrammatic notation of the later sections).

The $e_8$ series (which, suitably extended, includes all of the
exceptional Lie algebras) is the most intriguing. For $e_8$, the
smallest representation on which an $R$-matrix may be constructed
(and in fact the smallest representation of the Yangian $Y(e_8)$
\cite{drinf85}) is the $\g$-reducible representation $\g\oplus
\C$. Its $R$-matrix is constructed in \cite{chari91}, and Westbury
observes that this has a nice, uniform parametrization by Vogel's
plane \cite{vogel99}. (Note that such uniformity suggests an
extension of Deligne's conjecture \cite{delig96}, about the
uniformity of decomposition of $\g^{\otimes r}$, to Yangians.)
Although both conventional \cite{macfa04} and  diagrammatic
\cite{cvitabk} techniques for the adjoint representation of the
$e_8$ series (the latter as advocated in \cite{cohen96,westb03b})
are well-developed, we have not yet been able to extend them to
this reducible representation. Such a treatment of the $R$-matrix
remains, however, highly desirable, as a step towards explaining
the remarkable appearance of spectra associated with the algebras
of the $e_8$ series in the $q$-state Potts model $S$-matrix
\cite{dorey01,westb04}.

Westbury's observations also apply to trigonometric
($q$-dependent) $R$-matrices when $q$ is not a root of unity. As
far as we know, the centralizer algebras we study, which
$q$-deform to the Iwahori-Hecke algebra for the $su(n)$ and the
Birman-Wenzl-Murakami algebra for the other classical cases, have
not been constructed for exceptional $\g$ other than $g_2$
\cite{kuper96}.

\section{The classical Lie algebras}

Perhaps the two best-known, classic solutions of the YBE are those
of Yang \cite{yang67}, acting on the $n$-dimensional module of
$SU(n)$, \be \label{Yang}\check{R}^{ab}_{cd}(u) = 2
\delta^a_c\delta^b_d - u \delta^a_d
\delta^b_c=\left(2\,\I_n\otimes \I_n -
u\transp\right)^{ab}_{cd}\,, \ee and of the Zamolodchikovs
\cite{zamol79}, acting on the $n$-dimensional module of $SO(n)$,
\be\label{Zams} \check{R}^{ab}_{cd}(u) = 2 \delta^a_c\delta^b_d -
u \delta^a_d \delta^b_c + {2u \over n-2-u} \delta^a_b\delta^c_d\,.
\ee In the classic diagrammatic notation for these, which avoids a
proliferation of indices in calculations,  (\ref{Yang}) is
\be\label{Yang2} \Ru \;\;=\; 2 \;\soBA - u\; \soBB\ee and
(\ref{Zams}) is \be\label{Zams2} \Ru \;\;=\; 2 \;\soBA \;- u\;
\soBB\;+ {2u \over n-2-u}\;\soBC\;,\ee\\
 in which each Kronecker delta is written as a line connecting two indices.
Concatenation of symbols (by connecting lines, horizontally, from
right to left) is the correct way to multiply these (since
$\delta^a_b\delta^b_c=\delta^a_c$), so that the YBE becomes
\be\label{YBE3} \RUu \RDuv \RUv \;=\; \RDv \RUuv \RDu \;,\ee\vskip
-5pt\noindent in which internal lines represent summed indices and
external lines free indices. Checking that (\ref{Yang},\ref{Zams})
are indeed solutions is now a matter of checking the equivalence
of two $\C$-linear combinations of symbols, subject in the latter
case to the further condition that a loop takes value
$\delta^a_b\delta^b_a=n$.

\subsection{$su(n)$}

As already indicated, there is a Lie algebra and its
representation theory underlying each of these solutions. In the
first case, and denoting the $n$-dimensional module of $SU(n)$ by
$V$ and its conjugate by $\bar{V}$,  we re-write
(\ref{Yang},\ref{Yang2}) as \be \RuD \;\;= (2-u)P_+ + (2+u)P_-\ee
where \be\label{proj} P_\pm = {1\over
2}\left(\;\suBA\pm\suBB\;\right) \ee\vskip -5pt\noindent are
idempotents $P^2_\pm=P_\pm$ (and we henceforth distinguish $V$
from $\bar{V}$ by decorating each line with an arrow). In fact
these are the projectors onto the symmetric and antisymmetric
irreducible components of the tensor square $V\otimes V$, and we
thus have the spectral decomposition of the $R$-matrix, in form
(\ref{spectral}).

The $R$-matrix takes its values in the centralizer algebra
$End_{G}(V^{\otimes 2})$, the commutant of the action of the group
$G$ (and Lie algebra $\g$) on $V\otimes V$. The projectors
(\ref{proj}), therefore, or alternatively $\delta^a_c\delta^b_d$
and $\delta^a_d\delta^b_c$ and the symbols which represent them,
form a basis for $End_{su(n)}(V^{\otimes 2})=\C S_2$, the algebra
of the symmetric group $S_2$. Similarly the YBE is an equation of
expressions in $End_{su(n)}(V^{\otimes 3})=\C S_3$, or, in the
symbolic notation, $\C$-linear sums of the symbols \vskip -10pt\be
\suCA \;,\hsp \suCB\;,\hsp\suCC \;,\hsp \suCD\;,\hsp\suCE \hsp
{\rm and}\hsp \suCF\;. \ee\vskip-5pt\noindent This
mutually-centralizing action of $S_p$ and $SU(n)$ on $V^{\otimes
p}$ is the classic Schur-Weyl duality.

\subsection{$so(n)$}

We can rewrite (\ref{Zams},\ref{Zams2}) similarly as \be \Ru \;\;=
(2-u)P_+ + (2+u)P_- + (2+u){n-2+u\over n-2-u}P_0\,,\ee where \be
P_+  = {1\over 2}\left(\;\soBA + \soBB\;\right) - {1\over
n}\;\soBC \;,\hsp P_-  = {1\over 2}\left(\;\soBA -
\soBB\;\right)\;,\hsp P_0 = {1\over n} \;\soBC\ee are the
projectors onto the symmetric traceless, antisymmetric and singlet
components of the tensor square $V\otimes V$ of the defining,
$n$-dimensional representation of $so(n)$. To check that each
$P^2=P$, we need the algebraic relations among these symbols,
which are simply those of concatenation together with the loop
value $n$, or \vskip -15pt \beaa &
\soBB\soBB=\soBA\;,\hsp \soBA\soBB=\soBB\soBA=\soBB\;,\hspace{0.9in}& \\[0.3in]
&\hspace{-0.3in}
\soBA\soBC=\soBC\soBA=\soBB\soBC=\soBC\soBA=\soBC\;,\hsp
\soBC\soBC = n\;\soBC\;. \hspace{0.6in}& \\ \eeaa The dimension of
the module corresponding to the idempotent $P$ is computed in the
algebra by connecting the in- to the out- top index and the in- to
the out- bottom index, equivalent to taking the trace in the
tensor product by setting $a=c$ and $b=d$ and summing. This gives
values for $P_+,P_-$ and $P_0$ of $n(n+1)/2-1, n(n-1)/2$ and $1$
respectively.

This algebra, $End_{so(n)}(V^{\otimes 2})$, is  Brauer's algebra
$B_2(n)$ \cite{braue37,wenzl88}. The YBE is now valued in
$End_{so(n)}(V^{\otimes 3})=B_3(n)$, the 15-dimensional algebra
spanned by  \vskip -15pt \beaa & \soCA \spc\soCB \spc \soCC \spc
\soCD
\spc \soCE \spc \soCF & \\[0.35in] &
\soCG \spc \soCH \spc \soCI \spc \soCJ \spc \soCK \spc \soCL \spc
\soCM \spc \soCN \spc \soCO\;, &\\ \eeaa subject to the same rules
of concatenation and loop value $n$.

\subsection{$sp(2r)$}

There is another solution of the YBE \cite{berg78}, associated to
$sp(2r)$, which utilizes $B_2(-2r)$, although we shall instead
write it in a form which makes the role of the symplectic form
matrix explicit. It is \bea \Ru \;\; &=& (2-u)P_+ + (2+u)P_- +
(2-u){2r+2+u\over 2r+2-u}P_0\nonumber\\[0.1in] & = &  2 \;\spBA \;- u\;
\spBB\;+ {2u \over 2r+2-u}\;\spBC\;,\label{sp}\eea where\be P_+  =
{1\over 2}\left(\;\spBA + \spBB\;\right) \;,\hsp P_-  = {1\over
2}\left(\;\spBA - \spBB\;\right)+ {1\over 2r}\;\spBC \;,\hsp P_0 =
-{1\over 2r} \;\spBC\ee are the projectors onto the symmetric,
antisymmetric and symplectic-traceless, and singlet components of
$V\otimes V$. We use a solid arrow $\spAA$ to denote the
symplectic form matrix, so that $\spAA\spAA\;=-\,\soAA\,$,
$\;\spAA\;=-\,\spAAT$ and $\spAA\spAAT\;=\;\soAA$. If we denote an
element of $Sp(2r)$ by $\G$ then the defining relation $MJM^T=J$
for $M\in Sp(2r)$ is that \be\label{spdef} \G \spAA \GT \;=\;
\spAA\;.\ee The algebra $End_{sp(2r)}(V^{\otimes 2})$ is generated
by the three symbols in (\ref{sp}), with the invariance of the
third being due to (\ref{spdef}), \be \GG
\spBC\;=\;\spBC\;=\;\spBC\GG\;.\ee\vskip -5pt\noindent It is
simple to check that each of the three given projectors is indeed
idempotent. The YBE is valued in $End_{sp(2r)}(V^{\otimes 3})$,
spanned by \vskip -20pt\beaa & \spCA \spc\spCB \spc \spCC \spc
\spCD \spc
\spCE \spc \spCF & \\[0.3in] &
\spCG \spc \spCH \spc \spCI \spc \spCJ \spc \spCK \spc \spCL \spc
\spCM \spc \spCN \spc \spCO\;. & \eeaa

\subsection{Dimension of $End_{\g}(V^{\otimes p})$}

An alternative basis for the centralizer algebra
$End_{\g}(V^{\otimes p})$ (for semisimple $\g$) is given by the
set of projectors and intertwiners of $\g$-irreducible components
of $V^{\otimes p}=\bigoplus_i \C^{d_i} \otimes W_i$ (in which
$d_i$ is the multiplicity of $W_i$ in the decomposition). Thus
\be\label{dimcen} {\rm dim}\,End_{\g}(V^{\otimes p}) = \sum_i
d_i^2\,, \ee which we shall find useful in dealing with the
exceptional algebras, where a diagrammatic basis for
$End_{\g}(V^{\otimes 3})$ (as used in the previous subsections)
will be far from obvious. The central utility of such bases, which
is not achieved by using projectors and intertwiners, is to
facilitate calculations in $End_{\g}(V^{\otimes 3})$ using terms
from the different embeddings of $End_{\g}(V^{\otimes 2})$, as
required by the YBE.

\section{The $e_6$ series}

The defining property of $\g$ in the $e_6$ series, as subgroups
$G\subset SU(n)$ with $n=3m+3$, is the existence of a cubic,
symmetric invariant form $d_{abc}$, {\em i.e.\ }a map $V^{\otimes
3}\rightarrow\C$, $(u^a,v^b,w^c)\mapsto d_{abc}u^av^bw^c$ such
that, for $M\in G$, $d_{def} M^{da}M^{eb}M^{fc}=d_{abc}$ or, in
diagrammatic notation (and with $\G$ denoting $M$),
\be\label{e6D}\eCK\GGG \;= \eCK \hspace{0.3in}{\rm and}
\hspace{0.3in}
\eCK\suCB\;=\eCK\suCC\;=\eCK\suCD\;=\eCK\\\;.\ee\vskip-5pt\noindent
(This should not be confused with the cubic Casimir operator of
$su(n)$ corresponding to the cubic symmetric invariant in the {\em
adjoint} representation.) Following \cite{cvitabk}, this is
normalized so that \be \eBA\eBA \;=\;\eBA\;.\ee\vskip
-5pt\noindent Thus the symmetric component of $V\otimes V$
decomposes further, and the $R$-matrix (\ref{e6R}), with \be P_1=
{1\over 2}\left(\;\suBA\;+\;\suBB\;\right) - \eBA\;,\hsp P_2
={1\over 2}\left(\;\suBA\;-\;\suBB\;\right)\;,\hsp P_3 = \eBA \;,
\ee may be multiplied by $4-u$ to give \be\label{e6R2} \RuD \;=\;
4 \;\suBA \;- u\;\suBB + {(4m+8)u\over 2m-u} \;\eBA\;.\ee \vskip
-5pt\noindent (Note the re-scalings of $R$ and $u$ relative to the
$su(n)$ $R$-matrix (\ref{Yang},\ref{Yang2}).)

As discussed in section 3.2, one computes the putative trace of an
idempotent by connecting its in- and out- legs. For the
idempotents constructed using $d_{abc}$ this is an integer, and
thus the centralizer algebra has an action on a module, when
$n=3m+3$ for $m=1,2,4$ and 8 (although not only for these---for
the full story see \cite{cvitabk}). It is worth noting that the
centralizer algebras for all, including classical, $\g$ are
formally defined, and $R$-matrices in them exist, for {\em all}
$n$, not just integers: it is only the requirement that
idempotents have integer `trace' which further restricts $n$.

 The three-dimensional
centralizer $End_{\g}(V\otimes V)$ is generated by the three
symbols which appear in (\ref{e6R2}): note that the third symbol's
commuting with $E_6$ follows from (\ref{e6D}),\vskip -15pt
\be\eBA\GG \;=\; \eBAspec\;=\;\GG\eBA\;,\ee \vskip -5pt\noindent
in which the bar denotes complex conjugation, so that the defining
property of $U(n)$ is \be \G\GTbar\;=\;\suAA\;.\ee

Our object now is to demonstrate that (\ref{e6R}) satisfies the
YBE. It is clear that this will include terms of even orders in
$d$ up to six, and (for reasons which will become apparent below)
that there will be reduction relations among them. This will all
be rather involved, and so we move now to a more
mathematically-formal layout.

There are two primary identities satisfied by the invariant
$d_{abc}$, at fourth and third order respectively, and there are
no more at these or lower orders \cite{cvitabk}. The first is

\noindent {\bf Lemma 4.1} (Cvitanovi\'c):
$$\eBb\eBB \;= -{m\over 2m+4}\eBC + {1\over
2m+4}\left(\;\suBE\;+\;\suBF \;\right)\,.
$$
{\em Proof}  in \cite{cvitabk}, eqn.\ (18.9); follows from
irreducibility of components of $V\otimes \bar{V}$. \hfill$\Box$

All terms in the YBE are of rank six (where the rank, the number
of free indices, is the number of external legs of a diagram), and
to reduce the sixth-order terms in the YBE we will need\\[0.1in]
{\bf Corollary 4.2}:
 \vskip -0.3in
 \beaa
 \eCA\eCG\eCA\;-\;\eCG\eCA\eCG \;& = & {m\over 2m+4}\left(\;
 \eCG\suCB\eCG\;-\;\eCA\suCC\eCA\;\right)\\[0.1in] && \hspace{0.2in} - {1\over
 2m+4}\left(\; \eCA\;-\;\eCG \;\right)\,.
 \eeaa
 {\em Proof} is by applying Lemma 4.1 to the loops. \hfill$\Box$

With the sixth-order terms thus reduced, we now deal with the
fourth-order terms. To do so we begin with the other primary
identity satisfied by the invariant $d$,\\[0.1in]
{\bf Lemma 4.3} (Freudenthal):
$$
\left.\eDA\;\right\}{\rm symmetrized}\; = {4\over
3m+6}\left.\eDB\;\right\}{\rm symmetrized}.
$$\vskip -5pt\noindent
{\em Proof} in \cite{freud54}, eqn.\ (1.17). The cubic invariant
is the determinant of a $3\times 3$ hermitian matrix $X$ with
entries in the division algebra of order $m$ (and which thus form
$V$ of dimension $3m+3$). Freudenthal utilizes $d$ to define a
product $\times: V\otimes V\rightarrow \bar{V}$ (and its
conjugate) which obeys $(X\times X)\times (X\times X) = X{\rm det}
X$, expressed diagrammatically above. The relation appears
 in \cite{cvitabk} (sect.18.10) as the `Springer
relation' \cite{spring62}, and for $e_6$ specifically in
\cite{cvita76}, Fig.15(b). \hfill$\Box$

Once again we need relations of rank six rather than the
rank-five of Lemma 4.3, and so must construct the secondary identity\\[0.1in]
{\bf Corollary 4.4}:\vskip -25pt
$$
\eCA\eCG\;+\;\eCA\eCH\;+\;\eCA\suCC\eCA = {1\over
m+2}\left(\;\eCA\;+\;\eCB\;+\;\eCC\;+\;\eCJ\;\right)\,.
$$\vskip -5pt\noindent
{\em Proof} by appending another copy of $d$ to the diagrams of
Lemma 4.3, expanding the symmetrizers, and re-arranging.
\hfill$\Box$

These results are sufficient for us now to prove\\[0.1in]
{\bf Theorem 4.5}: $End_{e_6}(V^{\otimes 3})$ is the
20-dimensional algebra, with subalgebra $End_{su(27)}(V^{\otimes
3})=\C S_3$, spanned by
$$ \suCA  \spc \suCB \spc \suCC \spc
\suCD \spc \suCE \spc \suCF\;,$$ \vskip -0.1in
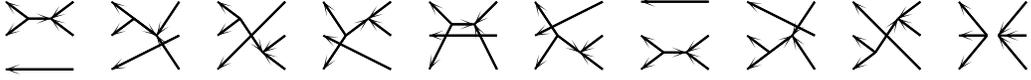
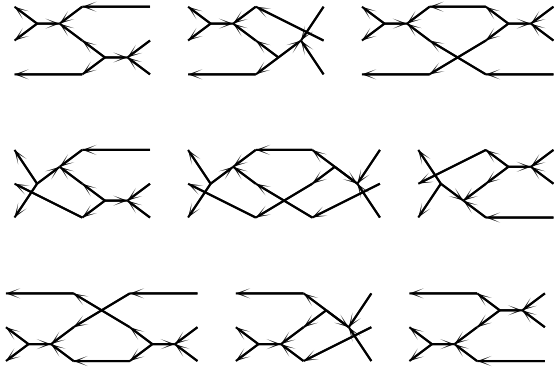
$$\eCA \spc \eCB \spc \eCC \spc \eCD
\spc \eCE \spc \eCF \spc \eCG \spc \eCH \spc \eCI \spc \eCJ$$\\
and
$$
\eCA\eCG\spc\eCA\eCH\spc\eCA\suCC\eCA
$$
\vskip -0.1in
$$
 \eCD\eCG \spc \eCD \suCC \eCB\spc\eCF\eCA
$$
\vskip -0.1in
$$
\eCG\suCB\eCG\spc \eCG\eCB \spc \eCG\eCA\;.
$$

\noindent{\em Proof.} The six symbols at zeroth order and the ten
at second order are trivially independent -- they cannot be
related by Lemma 4.3. At fourth order let us denote the nine
symbols by $e_{ij}$, $i,j=1,2,3$, where $i$ (respectively $j$)
indexes the left-hand, $V$ (resp.\ right-hand, $\bar{V}$) leg not
contracted on a common $d$. Corollary 4.4 (for $m=8$) then reduces
$e_{11}+e_{12}+e_{13}$ to terms of lower order. Permuting external
legs (and $\C$-conjugating where necessary) gives six such
reduction relations in total, reducing $e_{i1}+e_{i2}+e_{i3}$ and
 $e_{1i}+e_{2i}+e_{3i}$ for each $i=1,2,3$. Only five of these six
are independent, since $\sum_i \left( e_{i1}+e_{i2}+e_{i3} -
e_{1i} - e_{2i}- e_{3i}\right)=0$. There are therefore five
reduction relations among the nine symbols at fourth order,
leaving four independent generators. We thus have
dim$End_{\g}(V^{\otimes 3})=6+10+4=20$, matching that computed
from (\ref{dimcen}).\hfill$\Box$

\noindent{\bf Remark 4.6}. A set of four independent symbols among
the nine at fourth order is furnished by any set of four which
neither (i) includes three from any single row or column, nor (ii)
consists of two from one row and the other two from the excluded
column. An example is $\{e_{11},e_{12},e_{21},e_{22}\}$.

 \noindent{\bf Remark 4.7}. Theorem 4.5
does not apply to other $\g$ in the $e_6$ series, for which there
is a further reduction (which does not affect our YBE results).
For details, and an extended Young tableau method for the $e_6$
series, see
 ch.18 of \cite{cvitabk}.

For the YBE we will need some further fourth-order relations, for
which we begin with\\
{\bf Definition 4.8}: for any rank-six symbol $\scriptsize\obj$ we
define the transformations
$$ T_1:\spc \obj \spc\mapsto \spc\suCB\obj\suCC\spc,\hspace{0.5in}
T_2:\spc \obj \spc\mapsto \spc\suCC\obj\suCB\;,
$$\vskip 0pt\noindent
which facilitates\\[0.1in] \noindent{\bf Lemma 4.9}: the unique (up to
scaling) fourth-order term of rank six with eigenvalue $-1$ under
both $T_1$ and $T_2$ is
$$
\objj\;:=\;\eCG\suCB\eCG \;-\; \eCA\suCC\eCA
\;+\;\eCF\eCA\;-\;\eCG\eCB \;+\; \eCA\eCH \;-\;\eCD\eCG\,. $$ {\em
Proof} by direct calculation.\hfill$\Box$

Next is the key lemma in checking the YBE,\\[0.1in]
{\bf Lemma 4.10}: \beaa \hspace{0.2in}\eCG\suCB\eCG \;-\;
\eCA\suCC\eCA \;&=&\; {1\over 3}\;\objj \;+\;{2\over
3m+6}\left(\;\eCG\;-\;\eCA \;\right)
\\[0.1in] && \hspace{0.1in}+{1\over 3m+6}\left(\;  \eCF\;+\;\eCH\;-\;\eCD
\;-\;\eCB\;\right)\,.\eeaa {\em Proof.} This is a linear
combination of four of the six variants of Corollary 4.4. Using
again the basis $e_{11},...,e_{33}$ introduced in Theorem 4.5 for
the fourth-order terms, it is the reduction formula for ${1\over
3}\sum_i \left( e_{1i}+e_{i1}-e_{3i}-e_{i3}\right)$.\hfill$\Box$

It is such combinations of diagrams, and permutations (of the
external legs) thereof, which appear in the Yang-Baxter equation,
which we can now see is connected rather subtly, through the
secondary identities in Corollary 4.4 and Lemma 4.10, with
Freudenthal's primary relation, Lemma
4.3. Thus we can now prove\\[0.1in]
{\bf Theorem 4.11}: the $R$-matrix (\ref{e6R2}) solves the
YBE.\\[0.1in]
{\em Proof.} We first substitute (\ref{e6R2}) into the YBE
(\ref{YBE3}) and expand the left-hand- minus the right-hand-side.
The combination of sixth-order terms is precisely the
left-hand-side of Corollary 4.2, which thereby reduces the overall
expression to fourth-order. There are then three combinations of
fourth-order terms which appear,\vskip -20pt
$$
\eCG\suCB\eCG \;-\; \eCA\suCC\eCA \;,\hspace{0.25in} \eCD\eCG
\;-\;
 \eCA\eCH\hspace{0.15in}{\rm and} \hspace{0.15in}
\eCG\eCB\;-\;\eCF\eCA\,.$$\vskip 0pt\noindent To the first of
these we apply Lemma 4.10, and to the others, respectively, $T_1$
and $T_2$ of Lemma 4.10. The nice behaviour of $\scriptsize\objj$
under $T_1$ and $T_2$ ensures that its coefficient vanishes. What
remains is a linear combination of the zeroth- and second-order
symbols. That each of the coefficients vanishes was checked both
by hand and using Maple. \hfill$\Box$

\noindent{\bf Corollary 4.12}: the $V\otimes\bar{V}$ $R$-matrix
(\ref{e6Rcr}), with projectors \vskip -15pt $$ {\bf P}P_1=
{m+2\over m+4}\left(\;\suBC\;-{1\over
m+1}\;\suBD\;+2\eBB\;\right)\,, $$ \vskip -20pt$$
 {\bf P}P_2={2\over m+4}\left(\; \suBC\;+{1\over
3}\;\suBD\;-(m+2)\;\eBB \;\right)\hspace{0.2in}{\rm
and}\hspace{0.2in} {\bf P}P_3 = {1\over 3m+3} \;\suBD $$ \vskip
-10pt\noindent(from \cite{cvitabk}) and thus (rescaled)
$$
\RuDc \spc=\spc \suBC \;-4{m+2\over u-m}\;\eBB\; + {4\over
u-3m}\suBD \;,
$$\vskip -5pt\noindent
combines with the $V\otimes V$ $R$-matrix of Theorem 4.10 to solve
the YBE on $V\otimes V \otimes \bar{V}$. \\[0.1in]
{\em Proof}. We rely
here on the crossing-relation from factorized $S$-matrix theory
(see, for example, \cite{zamol80}), which in our case states that
$$ R_{V\bar{V}}(u) \propto Cross\left( R_{VV}(3m-u) \right)\,,
$$ where the operation $Cross$ simply rotates the symbolic
representation of $R$ anticlockwise through $90^\circ$. It is
simple to check that this holds, thereby implying that
(\ref{e6Rcr}) is indeed the correct $R$-matrix on
$V\otimes\bar{V}$. \hfill$\Box$

%\newpage
\section{The $e_7$ series}

The progression of ideas in this section is very similar to that
in the last. We begin by recalling that the defining property of
the $e_7$ series, realized as subgroups $G\subset Sp(2r)$, is the
existence of a symmetric, quartic invariant $d_{pqrs}$ in the
defining, $n$-dimensional module $V$ (where $n=2r=6m+8$). The
invariance, in diagrammatic notation, is \be\label{e7inv}
\eeDA\GGGG = \eeDA\,.\ee (In contrast to the last section, we use
a black disc to denote this quartic tensor, to avoid confusion
with the transposition diagram.) The idempotents have integer
trace here for dim$V=n=6m+8$ ($m=1,2,4,8$).

Once again the symmetric component of $V\otimes V$ now decomposes
further, modifying the projectors of section 3.3. The projectors
in the $R$-matrix (\ref{e7R}) (from \cite{cvitabk}, but here
rendered symbolically) are \vskip -10pt
$$
P_1={1\over 6(m+4)}\left\{\;
3(m+3)\left(\;\spBA\;+\;\spBB\;\right) -\;\eeBA \;\right\}\;,$$
\vskip -10pt$$ P_2 = {1\over 2}\left(\;\spBA\;-\;\spBB\;\right)
+{1\over 6m+8}\;\spBC\;,
$$\vskip -20pt
$$ P_3={1\over 6(m+4)}\left\{\;
3\left(\;\spBA\;+\;\spBB\;\right)+\;\eeBA\;\right\}\;,\hspace{0.5in}
P_4 =-{1\over 6m+8}\;\spBC\;.
$$ \vskip -5pt\noindent
The $R$-matrix (\ref{e7R}) is then, after re-scaling,
$$\Ru\spc=\; (2m+4-u)\;\spBA\; +u(u-m-1)\;\spBB\; + {u(2+u)\over
2m+2-u} \;\spBC\;+{u\over 3}\;\eeBA\,,$$ and we see that the
four-dimensional $End_{\g}(V\otimes V)$ for the $e_7$ series is
generated by these four symbols. The only subtlety is in combining
the symplectic form with $d$ in the last symbol: this is done so
that
$$
\GG \eeBA\;=\;\eeBA\GG\,,$$ \vskip -10pt\noindent in which we have
used both (\ref{e7inv}) and (\ref{spdef}).

Our object is to demonstrate that (\ref{e7R}) satisfies the YBE,
and a similar story of tensor identities to that for the $e_6$
series now follows. As before, there are two primary identities,
this time both of second order.

\noindent{\bf Lemma 5.1} (Cvitanovi\'c):$$ \eeBA\eeBA \;=
6(m+2)\;\eeBA + 18(m+3)\left(\; \spBA\;+\;\spBB \;\right)\,.$$
{\em Proof}: in \cite{cvitabk}, ch.\ 20, and specifically for
$e_7$ in \cite{cvita76}, Fig.18(e).\hfill$\Box$

Note that this, together with the relations of section 3.3 and
$$
\EEBA\spBC\;=0\;,\hspace{0.8in} \EEBA\spBB \;=\;\EEBA\,,
$$ fixes the structure of $End_{\g}(V^{\otimes 2})$.

\noindent{\bf Lemma 5.2} (Brown):
$$
\EECV\;+\;\EECQ\;-\;\EECP\;-\;\EECW \spc= \hspace{2in}$$
$$\hspace{1in}3\left(\; \EECJ\;-\;\EECL
\;+\;\EECG\;-\;\EECD\;+\;\EECN\;-\;\EECO\;+\;\EECI\;-\;\EECH
\;\right)$$ \vskip -5pt\noindent {\em Proof}: \cite{brown69},
section 3. Analogously to Freudenthal's relation for the $e_6$
series, Brown uses the quartic invariant to define an invariant
map $V^{\otimes 3}\rightarrow V$, of which this is the key
property. The primitive quartic invariant naturally occurs as the
contraction of the symplectic (2-)form with the alternating
(6-)form. The relation appears diagrammatically in \cite{cvitabk}
and \cite{cvita76}, Fig.15(d).\hfill$\Box$

These identities are sufficient to prove\\[0.1in]
{\bf Theorem 5.3}: $End_{e_7}(V^{\otimes 3})$ is the
35-dimensional algebra, with Brauer subalgebra
$End_{sp(56)}(V^{\otimes 3})$, spanned by \beaa & \spCA \spc\spCB
\spc \spCC \spc \spCD \spc
\spCE \spc \spCF\,, & \\[0.25in] &
\spCG \spc \spCH \spc \spCI \spc \spCJ \spc \spCK \spc \spCL \spc
\spCM \spc \spCN \spc \spCO\;, & \\[0.25in]& \eeCA \spc \eeCB \spc \eeCC
\spc \eeCD \spc \eeCE & \\[0.15in] &\eeCF \spc \eeCG \spc \eeCH
\spc \eeCI \spc \eeCJ & \\[0.15in] & \eeCK \spc \eeCL \spc \eeCM
\spc \eeCN \spc \eeCO & \eeaa and \beaa & \eeCP \spc \eeCQ \spc
\eeCR \spc \eeCS \spc \eeCT \spc \eeCU\;\;,&
\\[0.15in] &\;\; \eeCV \spc \eeCW \spc \eeCX \;\;, \spc\spc \eeCY\;\;.
& \eeaa \vskip -5pt\noindent {\em Proof.} The fifteen symbols of
zeroth order and the fifteen of first order are independent; they
cannot be related by Lemmas 5.1, 5.2. The terms at second order
are, however, subject to Lemma 5.2 and its variants obtained by
permuting external legs. Simple combinatorics superficially yields
twenty-four of these, six rotations multiplied by the four
possibilities of transposing (or not) the upper-left and
upper-right pairs of legs (the only pairs not already related by
symmetry in Lemma 5.2). However, Maple informs us that only five
of the 24 variants are independent. Thus the ten symbols at second
order are reduced by five independent reduction relations to five
independent symbols, and we have dim$End_{\g}(V^{\otimes
3})=15+15+5=35$, matching the computation
(\ref{dimcen}).\hfill$\Box$

\noindent{\bf Remark 5.4}. In contrast to the analogous result for
$e_6$ (Remark 4.6), we do not here have a neat general
characterization of all possible choices for five independent
terms among the ten at second order. However, from the form of
Lemma 5.2 and its variants it is straightforward to argue that
either (i) any one of the six first-row symbols together with the
four others, or (ii) any four of the six first-row symbols
together with any one of the next three, is likely to furnish an
independent set. That this is indeed so was checked, for all such
choices, using Maple.

 Before proving the key reduction
relations for the YBE, we first note\\[0.1in]
{\bf Lemma 5.5}:\vskip -20pt
$$
\objh\spc:=\spc \EECR\;+\;\EECS\;-\;\EECP\;-\; \EECU \;-\; \EECQ
\;+\; \EECT
$$\vskip -5pt\noindent
is the unique second-order, rank-six term
which is invariant under $60^\circ$ rotations.\\
\noindent{\em Proof} by direct calculation.\hfill$\Box$

The relations essential for the YBE are then\\[0.1in]
{\bf Lemma 5.6}\vskip -20pt
$$
\EECR\;-\;\EECU \spc={1\over 3}\;\objh \;+\;
2\left(\;\EECA\;-\;\EECB\;-\;\EECN\;-\;\EECO \;\right)$$\vskip
-0.5in
$$
\hspace{1.2in}+\;\EECE\;+\;\EECF\;-\;\EECD\;-\;\EECG\;-\;
\EECL\;-\;\EECM\;-\;\EECJ\;-\;\EECK
$$\vskip -5pt\noindent
{\em Proof.} Write $R_\theta$ for the anticlockwise rotation of a
symbol by angle $\theta$. Then this is $(R_{60^\circ} +
R_{120^\circ} - R_{240^\circ} - R_{300^\circ})$ of Lemma
5.2.\hfill$\Box$

\noindent{\bf Lemma 5.7}\vskip -30pt
$$\hspace{-0.2in}\EECZ\;\;-\;\;\EECZz\;=27(m+3)\left(\;
\SPCG\;-\;\SPCH\;+\;\SPCB\;-\;\SPCC\;+\;\SPCO\;+\;\SPCL\;-\;\SPCN\;-\;\SPCM
 \;\right)
$$
$$
\hspace{0.6in}-3(2m+5)\left\{
2\left(\EECN\;+\;\EECO\;+\;\EECF\;-\;\EECD\;+\;\EECE\;-\;\EECG\right)\right.
$$
$$\hspace{1.6in}\left.\EECK\;+\;\EECM\;+\;\EECJ\;+\;
\EECL\;+\;\EECA\;-\;\EECB\;\right\} \;-\;(m+1)\;\objh
$$\vskip -5pt\noindent
{\em Proof.} First we contract Lemma 5.6, on its top two indices,
with the bottom two indices of $\eeBAbs$, and use Lemma 5.1. This
gives the antisymmetrization of Lemma 5.7 on its bottom two
indices. Requiring $60^\circ$ rotational symmetry then forces the
required result. \hfill$\Box$

With these secondary identities of the invariant tensor established,
we can now prove\\[0.1in]
{\bf Theorem 5.8}:
 the $R$-matrix (\ref{e7R}) solves the
YBE.\\[0.1in]
{\em Proof.} Again we first expand the left-hand- minus
right-hand-side of the YBE (\ref{YBE3}) with (\ref{e7R})
substituted, and then multiply by $\SPCAA$. The third-order terms
are reduced by Lemma 5.7. Some of the second-order terms are
reduced by Lemma 5.1 alone; the others are the differences
$$
\EECR\;-\;\EECU\;\;,\hspace{0.5in} \EECP\;-\;\EECT
\hspace{0.3in}{\rm and}\hspace{0.3in}\EECQ\;-\;\EECS\;\;,
$$
or the left-hand side of Lemma 5.6 and its rotations by
$\pm120^\circ$. On using these the $\objh$ terms vanish because of
their invariant behaviour under rotations, Lemma 5.5. What remains
is an expression in the zeroth- and first-order, 30-dimensional
subalgebra of the centralizer. That each of the coefficients
vanishes was checked using Maple.\hfill$\Box$

\section{Concluding remarks}

In constructing and verifying the rational $R$-matrices for the
$e_6$ and $e_7$ series of Lie algebras, we have had to construct
their centralizers on $V^{\otimes 3}$ as diagram algebras, and
establish explicit bases for them (Theorem 4.5 for $e_6$ and
Theorem 5.3 for $e_7$). The connection between the algebras'
defining invariant tensors (and the primary reduction relations
satisfied by these) and the Yang-Baxter equation only appears
through a number of elegantly symmetric secondary identities
(Corollary 4.2 and Lemmas 4.4,4.10 for $e_6$, and Lemmas 5.6, 5.7
for $e_7$).

As we mentioned earlier, our primary goal for future work remains
to understand the $e_8$ case and its possible connections with the
$q$-state Potts model \cite{dorey01}.

 \vskip 0.1in

{\bf Acknowledgments}. NJM would like to thank Stephen Donkin for
discussions and comments on the draft, Tony Sudbery and James
Woodward for discussions, and Bruce Westbury for email
communications. AT would like to thank EPSRC for a PhD
studentship.

\end{document}